\documentclass[10pt,sumlimits,namelimits]{article}
\usepackage{amsthm}
\usepackage{amsmath}
\usepackage{amsfonts}
\usepackage{amssymb}
\usepackage{mathrsfs}
\usepackage{stmaryrd}
\usepackage{dsfont}
\usepackage{enumerate}
\usepackage{cancel}
\usepackage[applemac]{inputenc}
\usepackage[colorlinks=true,urlcolor=blue]{hyperref}

\newcommand{\1}{\mathds{1}}
\newcommand{\2}{{2^\star}}

\newcommand{\N}{\mathbb{N}}
\newcommand{\Z}{\mathbb{Z}}
\newcommand{\D}{\mathbb{D}}

\newcommand{\R}{\mathbb{R}}
\newcommand{\C}{\mathbb{C}}
\newcommand{\vsS}{\mathbb{S}}

\newcommand{\vN}{{\cal N}}

\newcommand{\Dr}{\mathscr{D}}

\newcommand{\Lr}{\mathscr{L}}
\newcommand{\Mr}{\mathscr{M}}

\newcommand{\vphi}{\varphi}
\newcommand{\eps}{\varepsilon}

\newcommand{\dsp}{\displaystyle}
\newcommand{\ovl}{\overline}
\newcommand{\udl}{\underline}
\newcommand{\vlim}{\lim\limits}
\newcommand{\vlimsup}{\limsup\limits}
\newcommand{\vliminf}{\liminf\limits}
\newcommand{\vmax}{\max\limits}

\newcommand{\vsup}{\sup\limits}
\newcommand{\vinf}{\inf\limits}
\newcommand{\vint}{\int\limits}

\newcommand{\inj}{\hookrightarrow}
\newcommand{\tends}{\longrightarrow}
\newcommand{\weak}{\rightharpoonup}
\newcommand{\wt}{\widetilde}

\newcommand{\li}{\llbracket}
\newcommand{\ri}{\rrbracket}
\newcommand{\loc}{\mathrm{loc}}

\newcommand{\co}{\mathrm{c}}

\renewcommand{\d}{\mathrm{d}}
\newcommand{\g}{\mathrm{g}}

\newcommand{\w}{{\textsl w}}
\renewcommand{\le}{\leqslant}
\renewcommand{\ge}{\geqslant}
\renewcommand{\Re}{\mathrm{Re}}
\renewcommand{\Im}{\mathrm{Im}}

\newcommand{\bs}{\boldsymbol}
\newcommand{\f}{\bs{f}}
\newcommand{\vi}{\mathrm{i}}
\newcommand{\p}{\prime}

\newcommand{\eqdef}{\stackrel{\text{def}}{=}}

\DeclareMathOperator{\supp}{supp}

\DeclareMathOperator{\essinf}{ess\:inf}
\DeclareMathOperator{\Id}{Id}

\DeclareMathOperator{\curl}{curl}

\numberwithin{equation}{section}

\newtheorem{thm}{Theorem}[section]
\newtheorem{prop}[thm]{Proposition}
\newtheorem{cor}[thm]{Corollary}
\newtheorem{lem}[thm]{Lemma}

\theoremstyle{definition}
\newtheorem{rmk}[thm]{Remark}
\newtheorem{defi}[thm]{Definition}

\newtheorem{exa}[thm]{Example}
\newtheorem{nota}[thm]{Notation}
\newtheorem{ass}[thm]{Assumption}

\newenvironment{proof*}{\noindent{\bf Proof.}}{\qed}
\newenvironment{vproof}[1]{\noindent{\bf Proof #1}}{\qed}

\title{On a Stationary Schrödinger Equation with Periodic Magnetic Potential}
\author{\sc Pascal Bégout\footnote{Institut de Mathématiques de Toulouse \& TSE, Universit\'{e} Toulouse I Capitole, 1 Esplanade de l'Université, 31080 Toulouse, Cedex 06, France} and Ian Schindler$^{*,\dagger}$}

\date{}

\begin{document}

\maketitle

\begin{abstract}
We prove existence results for a stationary Schrödinger equation with periodic magnetic potential satisfying a local integrability condition on the whole space using a critical value function.
\end{abstract}

{\let\thefootnote\relax\footnotetext{$^*$E-mail: \htmladdnormallink{{\footnotesize\udl{\tt{Pascal.Begout@math.cnrs.fr}}}}
{mailto:Pascal.Begout@math.cnrs.fr}}}
{\let\thefootnote\relax\footnotetext{$^\dagger$E-mail: \htmladdnormallink{{\footnotesize\udl{\tt{Ian.Schindler@ut-capitole.fr}}}}
{mailto:Ian.Schindler@ut-capitole.fr}}}
{\let\thefootnote\relax\footnotetext{Pascal Bégout and Ian Schindler acknowledge funding from ANR under grant ANR-17-EUR-0010 (Investissements d’Avenir program)}}
{\let\thefootnote\relax}
{\let\thefootnote\relax}

\tableofcontents

\baselineskip .65cm

\section{Introduction and main result}
\label{intro}

We wish to investigate for which $\lambda>0$ there is a weak solution  to
the stationary Schrödinger equation with magnetic potential:
\begin{gather}
\label{pde}
	\begin{cases}
          (-\vi\,\nabla + A)^2 u+ V(x)u= \lambda\f(x,u(x)), \; \text{ in } \R^N,	\medskip \\
          u\in H^1_{A,V}(\R^N),
	\end{cases}
\end{gather}
where $N\ge2,$ $A:\R^N\tends\R^N$ is the magnetic potential, $B=\curl A$ is the
magnetic field, $V:\R^N \tends\R,$ $\f$ satisfy some suitable assumptions, $\lambda>0$ and
\begin{gather}
\label{H1AV}
H^1_{A,V}(\R^N) \eqdef \Big\{u\in L^2(\R^N); Vu^2\in L^1(\R^N) \text{ and } (\nabla+\vi A)u\in L^2(\R^N;\C^N)\Big\}.
\end{gather}
Here, $\vi^2=-1$ and in what follows, unless specified, all functions are
complex-valued $(H^1(\R^N)=H^1(\R^N;\C),$ $L^p(\R^N)=L^p(\R^N;\C),$
$\Dr(\R^N)=\Dr(\R^N;\C),$ etc). If $N=1$ then~\eqref{pde} is equivalent to the case $A=0.$ Indeed, Assume that $A,V\in L^1_\loc(\R;\R)$ and $\f(\:.\:,e^{\vi\theta}z)=e^{\vi\theta}\f(\:.\:,z),$ for any $(\theta,z)\in\R\times\C$ (which is the case in this paper). Set for any $x\in\R,$ $\vphi(x)=\int_0^xA(s)\d s.$ If $u\in H^1_{A,V}(\R)$ is a solution to \eqref{pde} then by the gauge transformation $u\longmapsto v=e^{\vi\vphi}u,$ a straighforward calculation gives that $v\in H^1(\R)$ is a solution to $-\Delta v+V(x)v=\lambda\f(x,v(x)),$ which is~\eqref{pde} with $A=0.$ We thus restrict our study to the case $N\ge2.$

\medskip
\noindent
We make assumptions that insure the functional associated with \eqref{pde} is invariant with respect to the transformations $\tau_y:u \longmapsto e^{\vi\vphi_y} u(\:.\: +y),$ where $\varphi_y$ is defined in \eqref{defphi} and $y \in \Z^N.$ In \cite{MR1912753}, the authors stated that this set of transformations was a group of dislocations as defined in~\cite{MR2294665} which is false. In Section~\ref{sod} we prove (directly) that the set $D$ of such transformations is a set of dislocations permitting us to use the profile decomposition theorem \cite[Theorem 3.1, p.62-63]{MR2294665}. In Devillanova and Tintarev~\cite[Appendix]{MR4114959} this was proved by embedding the set of dislocations into a group via multiplication of $\tau_y\tau_z$ by $e^{\vi \alpha},$ $\alpha \in \R$ in such a way that the composition agrees with $\tau_{y+z}.$

\medskip
\noindent
Arioli and Szulkin~\cite{MR2022133} treated a similar problem with more
general conditions on $V$ (the spectrum of the operator $(-\vi\,\nabla +
A)^2 + V(x)$ can be negative), but they assume the Rabinowitz condition on
the right hand side.  We make less restrictive assumptions on the right hand
side and introduce a parameter $\lambda$ and an unbounded interval $I_\gamma\subset(0,\infty)$ such that for almost every $\lambda\in I_\gamma$
there is a solution to \eqref{pde}. In \cite{MR4114959} a magnetic Schrödinger equation with bounded non-periodic magnetic field is studied.

\medskip
\noindent
In Section~\ref{secequnor} we show that if the magnetic potential $A \in L^{N}_\loc(\R^N;\R^N)$ and $V \in
L^\frac{N}2_\loc(\R^N;\R)$ then $H^1_{A,V}(\R^N) = H^1(\R^N).$ In Section~\ref{sod}, we introduce the set of invariant dislocations acting
on \eqref{pde} and prove necessary results to the dislocation theorem in
\cite{MR2294665}. In Section~\ref{coc} we prove a cocompactness result. In
Section~\ref{scvf} we introduce a related critical
value function the study of which allows us to obtain our main result.  In
Section~\ref{secapp} we give some examples of nonlinearities to which our
result applies.

\noindent
Throughout this paper, we use the following notation. We denote by $\ovl z$ the conjugate of the complex number $z$ and by $\Re(z)$ its real part. For $1\le p\le\infty,$ $p^\p$ denotes the conjugate of $p$ defined by $\frac1p+\frac1{p^\p}=1.$ By $\{Q_j\}_{j \ge 1}$ we will denote a countable covering of $\R^N\setminus\Z^N$ by open unit cubes, thus $\R^N = \bigcup_{j \ge 1}\ovl{Q_j},$ and $Q=(0,1)^N.$ All vectors spaces considered will be over the field $\R.$ For a Banach space $X$ (over $\R),$ we denote by $X^\star\eqdef\Lr(X;\R)$ its topological dual and by $\langle\: . \; , \: . \:\rangle_{X^\star,X}\in\R$ the $X^\star-X$ duality product and for a Hilbert space $H,$ its (real) scalar product will be denoted by $\langle\:.\:,\:.\:\rangle_H.$ In particular, for any $T\in L^{p^\p}(\Omega)$ and $\vphi\in L^p(\Omega)$ with $1\le p<\infty,$ $\langle T,\vphi\rangle_{L^{p^\p}(\Omega),L^p(\Omega)}=\Re\vint_\Omega T(x)\ovl{\vphi(x)}\d x.$ See Appendix~\ref{appendixB} for more details. If $u\in L^p(\R^N),$ with $1\le p\le\infty,$ and if $\Omega$ is an open subset of $\R^N,$ with some abuse of notation, expression $\|u\|_{L^p(\Omega)}$ will stand for $\|u_{|\Omega}\|_{L^p(\Omega)}.$ This convention also holds for the others functional spaces. The subscript ``$\co$'' on a functional space indicates that the functions have compact support. For instance, if $\Omega\subseteq\R^N$ is an open subset then $u\in L^p_\co(\Omega)$ means that $u\in L^p(\Omega),$ $\supp u\subset\Omega$ and $\supp u$ is a compact subset of $\R^N.$ For a Banach space $E,$ the notation $E_\w$ designates the space $E$ endowed with the weak topology $\sigma(E,E^\star)$ and $E^\star_{\w\star}$ the space $E^\star$ endowed with the weak$\star$ topology $\sigma(E^\star,E).$ We denote by $C$ auxiliary positive constants which may change from a line to another one, and sometimes, for positive parameters $a_1,\ldots,a_n,$ write $C(a_1,\ldots,a_n)$ to indicate that the constant $C$ continuously depends only on $a_1,\ldots,a_n$ (this convention also holds for constants which are not denoted by ``$C$''). Finally, we denote by $\2=\frac{2N}{N-2}$ the critical exponent of the embedding $H^1(\R^N)\inj L^\2(\R^N),$ with the convention that $\2=\infty,$ if $N\le2.$

\medskip
\noindent
We shall make the following assumptions on $A:\R^N\tends\R^N.$

\begin{ass}
\label{assA}
Let  $(e_1,\ldots,e_n)$ be the canonical basis of $\R^N.$  
\begin{enumerate}
\item
\label{assA1}
The magnetic potential $A:\R^N\tends\R^N$ satisfies,
\begin{gather}
\label{A}
\begin{cases}
A\in L^N_\loc(\R^N;\R^N),							&	\text{if } N\ge3,	\medskip \\
A\in L^{2+\eps}_\loc(\R^2;\R^2), \text{ for some } \eps>0,	&	\text{if } N=2.
\end{cases}
\end{gather}
and
\begin{gather}
\label{A1}
\alpha_A\eqdef
\begin{cases}
\vsup_{j\in\N}\|A\|_{L^N(Q_j)}<\infty,			&	\text{if } N\ge3,	\medskip \\
\vsup_{j\in\N}\|A\|_{L^{2+\eps}(Q_j)}<\infty,	&	\text{if } N=2.
\end{cases}
\end{gather}
If $N\ge3$ then there exists $\eps>0$ such that,
\begin{gather}
\label{Ae}
A\in L^{N+\eps}_\loc(\R^N;\R^N).
\end{gather}
\item
\label{assA2}
$A$ is a \udl{$\Z^N-$periodic magnetic potential}:
\begin{gather}
\label{curlA}
\forall j\in\li1,N\ri, \; \curl A(x+e_j)\stackrel{\Dr^\p(\R^N)}{=}\curl A(x),
\end{gather}
where $\curl A\in\Mr_N\left(\Dr^\p(\R^N;\R)\right)$ is the skew-symmetric,
matrix-valued distribution with $A_{ij}=\partial_iA_j-\partial_jA_i.$
\end{enumerate}
\end{ass}

\begin{rmk}
\label{rmkhypA2}
It is easy to see that in Assumption~\ref{assA}, \eqref{curlA} is equivalent
to the condition: for any $y\in\Z^N,$
$\curl A(x+y)\stackrel{\Dr^\p(\R^N)}{=}\curl A(x).$ By Lemma~1.1 in
Leinfelder~\cite{MR695945}, \eqref{curlA} is also equivalent to: for any
$y\in\Z^N,$ there exists $\vphi_y\in W^{1,N+\eps}_\loc(\R^N;\R)$
such that for almost every
$x\in\R^N,$ $A(x+y)=A(x)+\nabla\vphi_y(x).$
\end{rmk}

\begin{ass}
\label{assf}
We will use the following
assumptions on $V$ and $f.$ Let $f:\R^N \times [0,\infty) \tends \R$ be such that $f(x,t)$ is measurable in $x$ and continuous in $t$ and let $F(x,t) \eqdef \dsp\int_0^t f(x,s) \d s,$ for almost every $x\in\R^N$ and any $t\ge0.$ We extend $f$ to the complex plane by setting for almost every $x\in\R^N$ and any $z\in\C\setminus\{0\},$ $\f(x,z)=f(x,|z|)\frac{z}{|z|},$ and $\f(x,0)=f(x,0).$ Finally, we set for any measurable function $u:\R^N\tends\C$ and almost every $x\in\R^N,$ $g(u)(x)=\f(x,u(x))$ and,
\begin{gather}
\label{psi}
\forall u\in H^1(\R^N), \; \psi(u) = \int_{\R^N} F(x,|u(x)|) \d x.
\end{gather}
\begin{enumerate}
\item
\label{assf1}
For every $\eps > 0,$ there exist $p_\eps\in(2,\2)$ and $C_\eps>0$ such that for almost every $x\in\R^N$ and any $t\ge0,$
\begin{gather}
\label{f1}
|f(x,t)| \le \eps (t + t^{\2-1}) + C_\eps t^{p_\eps-1},
\end{gather}
if $N \ge 3$ and
\begin{gather}
\label{f2}
|f(x,t)| \le \eps t + C_\eps t^{p_\eps-1},
\end{gather}
if $N\le2.$
\begin{gather}
\label{psi+}
\exists u\in H^1(\R^N) \; \text{ such that } \; \psi(u)>0.
\end{gather}
\item
\label{assf2}
The function $f$ and the electric potential $V:\R^N\tends\R$ are $\Z^N$-periodic, that is for almost every $(x,y)\in\R^N\times\Z^N$ and any $t\ge0,$ $f(x+y,t)=f(x,t)$ and $V(x+y) = V(x).$
\item
\label{assf3}
We have,
\begin{gather}
\label{V}
V\in L^1_\loc(\R^N;\R) \; \text{ and } \; \nu\eqdef\underset{x \in \R^N}{\essinf\;} V(x)>0.
\end{gather}
\end{enumerate}
\end{ass}

\begin{rmk}
\label{rmkf}
If $N\le2$ then Assumption~\eqref{f2} is equivalent to the following: 
\begin{gather}
\label{rmkf1}
\lim_{t\to0}\frac{f(x,t)}t=0, \; \text{ uniformly in } x\in\R^N,
\end{gather}
and there exist $p\in(2,\2)$ and $C>0$ such that,
\begin{gather}
\label{rmkf2}
\text{ for a.e. } x\in\R^N, \; \forall t\ge0, \; |f(x,t)|\le C+Ct^{p-1}.
\end{gather}
If $N\ge3$ then Assumption~\eqref{f1} is equivalent to~\eqref{rmkf1}, \eqref{rmkf2} with $p=\2$ and,
\begin{gather*}
\lim_{t\to\infty}\frac{f(x,t)}{t^{\2-1}}=0,
\end{gather*}
uniformly in $x\in\R^N.$
\end{rmk}

\begin{ass}
\label{assV}
The electric potential $V:\R^N\tends\R$ satisfies,
\begin{gather}
\label{V1}
\begin{cases}
V\in L^\frac N2_\loc(\R^N;\R) \text{ and } \alpha_V\eqdef\vsup_{j\in\N}\|V\|_{L^\frac N2(Q_j)}<\infty,	&	\text{if } N\ge3,
\medskip \\
V\in L^{1+\eps}_\loc(\R^2;\R) \text{ and } \alpha_V\eqdef\vsup_{j\in\N}\|V\|_{L^{1+\eps}(Q_j)}<\infty,
		\text{ for some } \eps>0,	&	\text{if } N=2.
\end{cases}
\end{gather}
\end{ass}

\begin{rmk}
\label{rmkV1}
Note that if $V$ has the local integrability~\eqref{V1} and if furthermore $V$ is $\Z^N$-periodic then we necessarily have $\alpha_V<\infty$ \big(since $\alpha_V=\|V\|_{L^1(Q)}\big).$
\end{rmk}

\begin{nota}
\label{notadaulHAV}
Let $A$ and $V$ satisfying \eqref{A}--\eqref{A1} and \eqref{V}, respectively, and let $H^1_{A,V}(\R^N)$ be defined by~\eqref{H1AV}. We shall denote by $H^{-1}_{A,V}(\R^N)$ the topological dual of the space $H^1_{A,V}(\R^N).$ This dual space is identified with a real vector subspace of the space of distributions $\Dr^\p(\R^N)$ (see Theorem~\ref{thmHHV} below).
\end{nota}

\begin{defi}
\label{defsol}
Let $A$ and $V$ satisfying \eqref{A}--\eqref{A1} and \eqref{V}, respectively, and let $H^1_{A,V}(\R^N)$ be defined by~\eqref{H1AV}. We shall write that $u$ is a \textit{weak solution} of~\eqref{pde} if $u\in H^1_{A,V}(\R^N)$ and if $u$ satisfies~\eqref{pde} in $H^{-1}_{A,V}(\R^N).$
\end{defi}

\begin{rmk}
\label{rmkdefsol}
With respect to Definition~\ref{defsol} we note:
\begin{enumerate}
\item
\label{rmkdefsol1}
If $u\in H^1(\R^N)$ then $(-\vi\nabla+A)^2u\in H^{-1}(\R^N)$ and,
\begin{align}
\label{rmkdefsol11}
&(-\vi\,\nabla + A)^2 u=-\Delta u-\vi\nabla.(Au)-\vi A.\nabla u+|A|^2u, \; \text{ in } H^{-1}(\R^N),				\\
\label{rmkdefsol12}
&\langle\vi A.\nabla u,v\rangle_{H^{-1}(\R^N),H^1(\R^N)}=\langle\vi \nabla u,Av\rangle_{L^2(\R^N),L^2(\R^N)},	\\
\label{rmkdefsol13}
&\langle|A|^2u,v\rangle_{H^{-1}(\R^N),H^1(\R^N)}=\langle Au,Av\rangle_{L^2(\R^N),L^2(\R^N)},
\end{align}
for any $v\in H^1(\R^N).$ Indeed, if $u\in H^1(\R^N)$ then by
Lemma~\ref{lemGA} below, $Au\in L^2(\R^N;\C^N)$ so that $-\Delta u\in
H^{-1}(\R^N)$ and $\nabla.(Au)\in H^{-1}(\R^N).$ In addition, by Hölder's inequality, $A.\nabla u,|A|^2u\in L^1_\loc(\R^N)\inj\Dr^\p(\R^N)$ and for any $\vphi\in\Dr(\R^N),$
\begin{gather*}
\langle\vi A.\nabla u,\vphi\rangle_{\Dr^\p,\Dr}=\langle\vi\nabla u,A\vphi\rangle_{L^2,L^2}
\; \text{ and } \;
\langle|A|^2u,\vphi\rangle_{\Dr^\p,\Dr}=\langle Au,A\vphi\rangle_{L^2L^2}.
\end{gather*}
By density and estimates in Property~\ref{lemps1} of Lemma~\ref{lemps} below, it follows that $A.\nabla u\in H^{-1}(\R^N),$ $|A|^2u\in H^{-1}(\R^N)$ and \eqref{rmkdefsol11}--\eqref{rmkdefsol13} follow.
\item
\label{rmkdefsol2}
Let $u\in H^1_{A,V}(\R^N).$ Let $R>0.$ We have,
\begin{gather*}
\int_{B(0,R)}|Vu|\d x\le\int_{B(0,R)\cap\{|u|\le1\}}|V||u|\d x+\int_{\{|u|>1\}}|V||u|^2\d x<\infty,
\end{gather*}
since $V\in L^1_\loc(\R^N;\R)$ and $Vu^2\in L^1(\R^N).$ It follows that $Vu\in L^1_\loc(\R^N)\inj\Dr^\p(\R^N)$ and for any $\vphi\in\Dr(\R^N),$
\begin{align*}
&\langle Vu,\vphi\rangle_{\Dr^\p,\Dr}=\Re\int_{\R^N}Vu\ovl\vphi\d x,		\\
&\left|\langle Vu,\vphi\rangle_{\Dr^\p,\Dr}\right|\le\|\sqrt Vu\|_{L^2}\|\sqrt V\vphi\|_{L^2}\le\|u\|_{H^1_{A,V}}\|\vphi\|_{H^1_{A,V}},
\end{align*}
by the Cauchy-Schwarz inequality (see Definition~\ref{defH1AV} below for the
definition of $\|\:.\:\|_{H^1_{A,V}}).$ By the density of $\Dr(\R^N)$ in $H^1_{A,V}(\R^N)$ (Theorem~\ref{thmHHV} below), it follows that $Vu\in H^{-1}_{A,V}(\R^N)$ and for any $v\in H^1_{A,V}(\R^N),$
\begin{gather}
\label{rmkdefsol21}
\langle Vu,v\rangle_{H^{-1}_{A,V}(\R^N),H^1_{A,V}(\R^N)}=\Re\int_{\R^N}Vu\,\ovl v\,\d x.
\end{gather}
Finally, by Proposition~\ref{propg} below, $g(u)\in H^{-1}(\R^N).$ In conclusion, since $H^{-1}(\R^N)\inj H^{-1}_{A,V}(\R^N)$ (Theorem~\ref{thmHHV} below), it follows from \eqref{rmkdefsol11} and \eqref{rmkdefsol21} that
\begin{gather*}
(-\vi\,\nabla+A)^2 u\in H^{-1}_{A,V}(\R^N),	\; Vu\in H^{-1}_{A,V}(\R^N) \; \text{ and } \; g(u)\in H^{-1}_{A,V}(\R^N).
\end{gather*}
Thus Definition~\ref{defsol} makes sense.
\end{enumerate}
\end{rmk}

\noindent 
Our  main result follows.
\begin{thm}
\label{thmmain}
Let $N\ge2$ and let Assumptions~$\ref{assA}$ and $\ref{assf}$ be satisfied. Let
$H^1_{A,V}(\R^N)$ be defined by~\eqref{H1AV}. Then for almost every
$\lambda>0$ sufficiently large, there exists, at least one non zero weak solution to,
\begin{gather}
\label{nlsp}
	\begin{cases}
          -\Delta_Au+V(x)u= \lambda g(u) \; \text{ in } \R^N,	\medskip \\
          u\in H^1_{A,V}(\R^N),
	\end{cases}
\end{gather}
where $-\Delta_Au=(-\vi\,\nabla + A)^2u.$
\end{thm}

\section{The space $\bs{H^1_{A,V}(\R^N)}$ and an equivalent definition of $\bs{H^1(\R^N)}$}
\label{secequnor}

In this section, we study the $H^1_{A,V}(\R^N),$ including the one-dimensional case $N=1$ because we believe that it is of interest for itself. For $N=1,$ the corresponding assumptions to~\eqref{A1} and \eqref{V1} are
\begin{gather}
\label{AN1}
A\in L^2_\loc(\R;\R) \text{ and } \alpha_A\eqdef\vsup_{j\in\N}\|A\|_{L^2(Q_j)} < \infty, \\
\label{VN1}
V\in L^1_\loc(\R;\R) \text{ and } \alpha_V\eqdef\vsup_{j\in\N}\|V\|_{L^1(Q_j)} < \infty,
\end{gather}
respectively.

\begin{defi}
\label{defH1AV}
Let $N\ge1$ and let $A\in L^2_\loc(\R^N;\R^N)$ and $V\in L^1_\loc(\R^N;\R)$ satisfy \eqref{A} and \eqref{V}, respectively. We recall that $H^1_{A,V}(\R^N)$ is defined by,
\begin{gather}
\label{H1AV1}
H^1_{A,V}(\R^N) = \Big\{u\in L^2(\R^N); Vu^2\in L^1(\R^N) \text{ and } \nabla_Au\in L^2(\R^N;\C^N)\Big\}.
\end{gather}
where $\nabla_Au=(\nabla+\vi A)u.$ We endow $H^1_{A,V}(\R^N)$ with the following scalar product and its corresponding norm,
\begin{gather*}
\forall u,v\in H^1_{A,V}(\R^N), \;
\langle u,v\rangle_{H^1_{A,V}(\R^N)}=\Re\vint_{\R^N}Vu\,\ovl v\,\d x+\Re\vint_{\R^N}\nabla_Au.\ovl{\nabla_Av}\d x, \\
\forall u\in H^1_{A,V}(\R^N), \;
\|u\|_{H^1_{A,V}(\R^N)}^2=\langle u,u\rangle_{H^1_{A,V}(\R^N)}=\vint_{\R^N}V|u|^2\d x+\|\nabla_Au\|_{L^2(\R^N)}^2,
\end{gather*}
making this space a real pre-Hilbert space. Indeed, it follows from~\eqref{V} that $\langle\;.\;,\;.\;\rangle_{H^1_{A,V}(\R^N)}$ is a bilinear symmetric positive definite form on $H^1_{A,V}(\R^N)\times H^1_{A,V}(\R^N).$
\end{defi}

\begin{rmk}
\label{rmkdefH1A}
Below are some comments about the definition of $H^1_{A,V}(\R^N).$
\begin{enumerate}
\item
\label{rmkdefH1A1}
If $u\in H^1_{A,V}(\R^N)$ then $\nabla_Au\eqdef(\nabla+\vi A)u\in L^2(\R^N;\C^N)$ but, \textit{a priori,} we do not assume that $\nabla u$ or $Au$ belong separately in $L^2(\R^N).$
\item
\label{rmkdefH1A2}
Frequently, in the literature (see for instance Sections~7.19--7.22,
p.191--195, of~Lieb and Loss~\cite{MR1817225}), it is assumed that $A\in
L^2_\loc(\R^N;\R^N)$ rather than
$A\in L^N_\loc(\R^N;\R^N)$  and $V\equiv1.$ With these assumptions it can
be shown that $H^1_{A,1}(\R^N)$ is a Hilbert space having $\Dr(\R^N)$ as a dense
subset. Moreover, if $u\in H^1_{A,1}(\R^N)$ then $|u|\in H^1(\R^N)$ and the
so-called \textit{diamagnetic inequality} \eqref{thmdi1} below holds.
However if $A \not \in L^N_\loc(\R^N;\R^N)$ then $H^1(\R^N)\not\subset H^1_{A,1}(\R^N)$ and
$H^1_{A,1}(\R^N)\not\subset H^1(\R^N).$ We show that if $A\in
L^N_\loc(\R^N;\R^N)$ then $H^1_{A,1}(\R^N)=H^1(\R^N)$ (see Theorem~\ref{thmHH} below).
\item
\label{rmkdefH1A3}
  Arioli and Szulkin showed (Lemma~2.3 in~\cite{MR2022133}) that if $N\ge2$
  and $A\in L^N_\loc(\R^N;\R^N)$ $(A\in L^{2+\eps}_\loc(\R^N;\R^N),$ if
  $N=2)$ then  $H^1_{A,1}(\Omega)=H^1(\Omega)$ with equivalent norms for
  open bounded subsets $\Omega$ of $\R^N$ with smooth boundaries. We extend
  their result to the case $\Omega=\R^N$ for any $N\ge1,$ under assumptions~\eqref{A}--\eqref{A1} (Theorem~\ref{thmHH} below).
\end{enumerate}
\end{rmk}

\begin{thm}
\label{thmHHV}
Let $A\in L^2_\loc(\R^N;\R^N)$ and $V\in L^1_\loc(\R^N;\R)$ satisfy \eqref{A}--\eqref{A1} $(\eqref{AN1},$ if $N=1)$ and \eqref{V}, respectively, and let $H^1_{A,V}(\R^N)$ be defined by~\eqref{H1AV1}. Then,
\begin{gather}
\label{thmHHV1}
H^1_{A,V}(\R^N)=\Big\{u\in H^1(\R^N); Vu^2\in L^1(\R^N) \text{ and } Au\in L^2(\R^N;\C^N)\Big\},	\\
\label{thmHHV2}
H^1_{A,V}(\R^N) \text{ is a separable Hilbert space,}										\\
\label{thmHHV3}
\Dr(\R^N)\inj H^1_{A,V}(\R^N)\inj H^1(\R^N),
\end{gather}
with both  dense embeddings. In particular, each term in the integrals of $\langle\:.\:,\:.\:\rangle_{H^1_{A,V}(\R^N)}$ belongs to $L^1(\R^N).$ In addition,
\begin{gather}
\label{thmHHV4}
H^{-1}(\R^N)\inj H^{-1}_{A,V}(\R^N)\inj\Dr^\p(\R^N),
\end{gather}
where $H^{-1}_{A,V}(\R^N)\eqdef\left(H^1_{A,V}(\R^N)\right)^\star$ and  both dense embeddings. 
\end{thm}

\begin{rmk}
\label{rmkthmHHV}
By~\eqref{thmHHV3}, for any $u\in H^1_{A,V}(\R^N),$ there exists $(\vphi_n)_{n\in\N}\subset\Dr(\R^N)$ such that $\vphi_n\xrightarrow[n\to\infty]{H^1_{A,V}}u.$ As with the classical proofs of density  we have for any $n\in\N,$ $\|\vphi_n\|_{L^p(\R^N)}\le\|u\|_{L^p(\R^N)},$ and $p\in[1,\infty],$ if $u\in L^p(\R^N).$ See Lemma~\ref{lemHHV2}.
\end{rmk}

\begin{thm}
\label{thmHH}
Let $A\in L^2_\loc(\R^N;\R^N)$ satisfy \eqref{A}--\eqref{A1} $(\eqref{AN1},$ if $N=1)$ and let $V\in L^1_\loc(\R^N;\R)$ satisfy \eqref{V} and \eqref{V1} $(\eqref{V}$ and $\eqref{VN1},$ if $N=1).$ Then,
\begin{gather*}
H^1_{A,V}(\R^N)=H^1(\R^N),
\end{gather*}
with equivalent norms.
\end{thm}

\begin{rmk}
\label{rmkH1AV}
To find  examples such that $H^1_{A,V}(\R^N)\subsetneq H^1(\R^N)$, note that
by~\eqref{lemGA2} below, assuming \eqref{A}--\eqref{A1}, for any $u\in H^1(\R^N),$ $u\in H^1_{A,V}(\R^N)$
if, and only if, $Vu^2\in L^1(\R^N).$ So we look for a $V$ that does not
satisfy~\eqref{V1}. If $N=1$ we must have $\alpha_V=\infty.$ In other words,
$V$ cannot be $\Z^N$-periodic. Below, for each $N,$ we give an example of a
$V$ and  a  $u\in H^1(\R^N)$ such that $u\not\in H^1_{A,V}(\R^N).$ For
$N\ge2,$ $u$ is a positive continuous function over $\R^N\setminus\{0\}$ such that for $|x|>10,$ $u(x)=|x|^{-N}.$ We give its definition for $x$ near $0$ below.
\begin{enumerate}
\item
\label{rmkH1AV1}
For $N=1.$ Define for any $x\in\R,$ $V(x)=x^4+1$ and $u(x)=(x^2+1)^{-1}.$ Then $V$ satisfies~\eqref{V}, $u\in H^1(\R^N;\R)$ but $Vu^2\not\in L^1(\R^N;\R).$
\item
\label{rmkH1AV2}
For $N=2.$ For $|x|<e^{-e},$ let $u(x)=\ln|\ln|x||$ and let $V$ be $\Z^N$-periodic such that for any $x\in Q,$
\begin{gather*}
V(x)=\frac1{|x|^2|\ln|x||(\ln|\ln|x||)^2}\1_{\left\{0<|x|<e^{-e}\right\}}(x)+\1_{\left\{|x|\ge e^{-e}\right\}\cap Q}(x).
\end{gather*}
Then $V$ satisfies~\eqref{V} but for any $p\in(1,\infty],$ $V\not\in L^p_\loc(\R^2;\R),$ so that $V$ does not verify \eqref{V1}. In addition, $u\in H^1(\R^N;\R)$ but $Vu^2\not\in L^1(\R^N;\R).$
\item
\label{rmkH1AV3}
For $N\ge3.$ Let for $|x|<1,$ $u(x)=|x|^{-\frac{N-2}4}$ and let $V$ be $\Z^N$-periodic such that for any $x\in Q,$ $V(x)=|x|^{-\frac{N+2}2}.$ Then $V$ satisfies~\eqref{V} but $V\not\in L^\frac{N}2_\loc(\R^N;\R),$ so that $V$ does not verify \eqref{V1}. In addition, $u\in H^1(\R^N;\R)$ but $Vu^2\not\in L^1(\R^N;\R).$
\end{enumerate}
\end{rmk}

\noindent
We recall that $-\Delta_A=(-\vi\,\nabla + A)^2$ and $\nabla_A=\nabla+\vi A.$
\begin{thm}
\label{thmLA}
Let $A\in L^2_\loc(\R^N;\R^N)$ satisfy \eqref{A}--\eqref{A1} $(\eqref{AN1},$ if $N=1).$ If $u\in H^1(\R^N)$ then $-\Delta_Au\in H^{-1}(\R^N)$ and for any $v\in H^1(\R^N),$
\begin{gather}
\label{thmLA1}
\langle-\Delta_Au,v\rangle_{H^{-1}(\R^N),H^1(\R^N)}=\langle\nabla_Au,\nabla_Av\rangle_{L^2(\R^N),L^2(\R^N)}.
\end{gather}
If in addition $V$ satisfies \eqref{V} and if $H^1_{A,V}(\R^N)$ is defined by~\eqref{H1AV1} then for any $T\in H^{-1}(\R^N),$ $T\in H^{-1}_{A,V}(\R^N)$ and
\begin{gather}
\label{thmLA2}
\langle T,u\rangle_{H^{-1}_{A,V}(\R^N),H^1_{A,V}(\R^N)}=\langle T,u\rangle_{H^{-1}(\R^N),H^1(\R^N)}.
\end{gather}
for any $u\in H^1_{A,V}(\R^N),$ 
\end{thm}

\begin{rmk}
\label{rmkthmmain}
Let $\lambda>0$ and let $u$ be a solution to \eqref{nlsp}. By Definition~\ref{defsol}, we may take the $H^{-1}_{A,V}-H^1_{A,V}$ duality product of \eqref{nlsp} with $u.$ We have by \eqref{thmLA2}, \eqref{thmLA1} and \eqref{rmkdefsol21},
\begin{gather}
\label{rmkthmmain1}
\|u\|_{H^1_{A,V}(\R^N)}^2=\lambda\langle g(u),u\rangle_{H^{-1}_{A,V}(\R^N),H^1_{A,V}(\R^N)}\ge0.
\end{gather}
It follows that if $g\equiv0$ then necessarily $u\equiv0.$ Note that if $\psi,$ defined by~\eqref{psi}, satisfies~\eqref{psi+} then $g\not\equiv0.$
\medskip
\end{rmk}

\begin{rmk}
\label{rmkDeltalin}
Let $A\in L^2_\loc(\R^N;\R^N)$ satisfy \eqref{A}--\eqref{A1} $(\eqref{AN1},$ if $N=1).$ It follows from Theorem~\ref{thmLA} and \eqref{lemGA3} below that $-\Delta_A\in\Lr\big(H^1(\R^N);H^{-1}(\R^N)\big).$ If, in addition, $V\in L^1_\loc(\R^N;\R)$ satisfies \eqref{V} then by Theorems~\ref{thmHHV} and \ref{thmLA}, $-{\Delta_A}_{|H^1_{A,V}(\R^N)}\in\Lr\big(H^1_{A,V}(\R^N);H^{-1}_{A,V}(\R^N)\big).$
\medskip
\end{rmk}

\noindent
We split the proof of Theorem~\ref{thmHHV} in several lemmas. We begin by
recalling the \textit{diamagnetic inequality} for functions belonging in
$H^1_{A,V}(\R^N).$ Its proof is well-known (Lieb and Loss~\cite{MR1817225},
Theorem~7.21, p.193). For the sake of completeness, we sketch the proof.

\begin{thm}[\textbf{Diamagnetic inequality, \cite{MR1817225}}]
\label{thmdi}
Let $A\in L^2_\loc(\R^N;\R^N).$ Let $u\in L^2(\R^N)$ be such that $(\nabla+\vi A)u\in L^2(\R^N;\C^N).$ Then $\nabla u,A u\in L^1_\loc(\R^N;\C^N),$ $|u|\in H^1(\R^N),$ $(\nabla+\vi A)u\stackrel{\text{a.e.}}{=}\nabla u+\vi Au$ and
\begin{gather}
\label{thmdi1}
|\nabla|u|\,|\stackrel{\text{a.e.}}{\le}|\nabla u+\vi Au|.
\end{gather}
\end{thm}

\begin{lem}[\textbf{\cite{MR1817225}}]
\label{lemH1}
If $u\in H^1(\R^N)$ then $|u|\in H^1(\R^N;\R)$ and $|\nabla|u|\,|\stackrel{\text{a.e.}}{\le}|\nabla u|.$
\end{lem}

\begin{vproof}{of Theorem~\ref{thmdi} and Lemma~\ref{lemH1}.}
We recall that if $u\in W^{1,p}_\loc(\R^N),$ for some $1\le p\le\infty,$
then $|u|\in W^{1,p}_\loc(\R^N;\R)$ and
$\nabla|u|\stackrel{\text{a.e.}}{=}\Re\left(\frac{\ovl u}{|u|}\nabla
  u\right)\footnote{$\nabla|u|=0,$ almost everywhere where $u=0.$}$
(Theorem~6.17, p.152, in~Lieb and Loss~\cite{MR1817225}). This proves
Lemma~\ref{lemH1}. Now, let $u\in L^2(\R^N)$ be such that $(\nabla+\vi
A)u\in L^2(\R^N;\C^N).$ Then $\nabla u\in H^{-1}(\R^N;\C^N)$ and, by the Cauchy-Schwarz inequality, $Au\in L^1_\loc(\R^N;\C^N).$ This implies that $\nabla u\in L^1_\loc(\R^N;\C^N).$ We then infer, $u\in W^{1,1}_\loc(\R^N)$ and $(\nabla+\vi A)u\stackrel{\text{a.e.}}{=}\nabla u+\vi Au.$ And since $\Re\left(\frac{\ovl u}{|u|}(\nabla u+\vi Au)\right)\stackrel{\text{a.e.}}{=}\Re\left(\frac{\ovl u}{|u|}\nabla u\right)\stackrel{\text{a.e.}}{=}\nabla|u|,$ one obtains that $|u|\in H^1(\R^N)$ and \eqref{thmdi1}.
\medskip
\end{vproof}

\begin{lem}
\label{lemA}
Let $u\in L^2(\R^N)$ be such that $|u|\in H^1(\R^N;\R).$
\begin{enumerate}
\item
\label{lemiA1}
If $A\in L^2_\loc(\R^N;\R^N)$ satisfies \eqref{A}--\eqref{A1} $(\eqref{AN1},$ if $N=1)$ then $Au\in L^2(\R^N;\C^N)$ and,
\begin{gather}
\label{lemA1}
\|Au\|_{L^2(\R^N)}\le C \alpha_A \|\,|u|\,\|_{H^1(\R^N)},
\end{gather}
where $C=C(N)$ $(C=C(\eps),$ if $N=2).$
\item
\label{lemiV1}
If $V$ satisfies \eqref{V1} $(\eqref{VN1},$ if $N=1)$ then $Vu^2\in L^1(\R^N)$ and,
\begin{gather}
\label{lemV1}
\|Vu^2\|_{L^1(\R^N)}\le C\alpha_V\|\,|u|\,\|_{H^1(\R^N)}^2,
\end{gather}
where $C=C(N)$ $(C=C(\eps),$ if $N=2).$
\end{enumerate}
\end{lem}

\begin{proof*}
Let $u\in L^2(\R^N)$ be such that $|u|\in H^1(\R^N;\R).$ We start by proving Property~\ref{lemiA1} with $N\ge3.$ By the Sobolev embedding $H^1(Q_j)\inj L^\2(Q_j),$ there exists $C=C(N,|Q_j|)$ such that for any $j\in\N,$ $\|u\|_{L^\2(Q_j)}\le C\|\,|u|\,\|_{H^1(Q_j)}.$ Actually, $C$ only depends on $N$ since for any $j\in\N,$ $|Q_j|=1.$ It follows from Hölder's inequality that,
\begin{eqnarray*}
\int_{\R^N}|Au|^2\d x &  =	&	\sum_{j\in\N}\int_{Q_j} |Au|^2 \d x					\\
				& \le	&	\sum_{j\in\N} \|A\|_{L^N(Q_j)}^2 \|u\|_{L^\2(Q_j)}^2		\\
				& \le & 	C^2\alpha_A^2\sum_{j\in\N}\|\,|u|\,\|_{H^1(Q_j)}^2		\\
				&  =	&	C^2\alpha_A^2\|\,|u|\,\|_{H^1(\R^N)}^2.
\end{eqnarray*}
If $N=2$ then the second line is replaced with $\sum_{j\in\N} \|A\|_{L^{2+\eps}(Q_j)}^2\|u\|_{L^\frac{2(2+\eps)}\eps(Q_j)}^2$ and we use the embedding $H^1(Q_j)\inj L^\frac{2(2+\eps)}\eps(Q_j),$ while if $N=1$ then the second line is replaced with
$\sum_{j\in\N} \|A\|_{L^2(Q_j)}^2 \|u\|_{L^\infty(Q_j)}^2$ and we use the embedding $H^1(Q_j)\inj L^\infty(Q_j).$ Hence \ref{lemiA1}. Property~\ref{lemiV1} follows in the same way: replace $A$ with $\sqrt{|V|}$ in the above estimates.
\medskip
\end{proof*}

\begin{lem}
\label{lemGA}
Let $A\in L^2_\loc(\R^N;\R^N)$ satisfy \eqref{A}--\eqref{A1} $(\eqref{AN1},$ if $N=1).$ Then,
\begin{gather}
\label{lemGA1}
(\nabla+\vi A)u\in L^2(\R^N;\C^N) \iff u\in H^1(\R^N),		\\
\label{lemGA2}
\nabla u\in L^2(\R^N;\C^N) \implies A u\in L^2(\R^N;\C^N).
\end{gather}
Finally, if $u\in H^1(\R^N)$ then  $(\nabla+\vi A)u=\nabla u+\vi Au,$ in $L^2(\R^N;\C^N)$ and
\begin{gather}
\label{lemGA3}
\|\nabla u+\vi Au\|_{L^2(\R^N)}\le C\|u\|_{H^1(\R^N)},
\end{gather}
where $C=C(\alpha_A,N)$ $(C=C(\alpha_A,N,\eps),$ if $N=2).$
\end{lem}

\begin{proof*}
Let $A$ satisfy \eqref{A}--\eqref{A1} (\eqref{AN1}, if $N=1).$ Let $u\in L^2(\R^N).$

$\bullet$
If $\nabla u\in L^2(\R^N;\C^N)$ then by Lemmas~\ref{lemH1} and \ref{lemA}, $Au\in L^2(\R^N;\C^N).$ Hence \eqref{lemGA2} and $\Longleftarrow$ in \eqref{lemGA1}.

$\bullet$
If $(\nabla+\vi A)u\in L^2(\R^N;\C^N)$ then by Theorem~\ref{thmdi} and Lemma~\ref{lemA}, $Au\in L^2(\R^N;\C^N)$ and $(\nabla+\vi A)u\stackrel{\text{a.e.}}{=}\nabla u+\vi Au.$ Hence, $\nabla u\in L^2(\R^N;\C^N)$ and $\implies$ in \eqref{lemGA1} is proved.

$\bullet$
By \eqref{lemA1} and Lemma~\ref{lemH1}, we have
\begin{gather*}
\|\nabla u+\vi Au\|_{L^2(\R^N)}\le\|\nabla u\|_{L^2(\R^N)}+\|Au\|_{L^2(\R^N)}\le(C\alpha_A+1)\|u\|_{H^1(\R^N)}.
\end{gather*}
Hence the result.
\medskip
\end{proof*}

\begin{lem}
\label{lemps}
Let $u,v\in H^1(\R^N).$
\begin{enumerate}
\item
\label{lemps1}
Let $A\in L^2_\loc(\R^N;\R^N)$ satisfy \eqref{A}--\eqref{A1} $(\eqref{AN1},$ if $N=1).$ Then $(Au).\nabla v\in L^1(\R^N),$ $|A|^2uv\in L^1(\R^N)$ and we have,
\begin{align*}
&	\int_{\R^N}|Au|\,|\nabla v|\d x\le C\alpha_A\|u\|_{H^1(\R^N)}\|v\|_{H^1(\R^N)},	\\
&	\int_{\R^N}|A|^2|uv|\d x\le C^2\alpha_A^2\|u\|_{H^1(\R^N)}\|v\|_{H^1(\R^N)}.
\end{align*}
where the constant $C$ is given by~\eqref{lemA1}.
\item
\label{lemps2}
Let $V$ satisfy \eqref{V1} $(\eqref{VN1},$ if $N=1).$ Then $Vuv\in L^1(\R^N)$ and we have,
\begin{gather*}
\int_{\R^N}|V||uv|\d x\le C\alpha_V\|u\|_{H^1(\R^N)}\|v\|_{H^1(\R^N)},
\end{gather*}
where the constant $C$ is given by~\eqref{lemV1}.
\end{enumerate}
\end{lem}

\begin{proof*}
The results come from Lemma~\ref{lemH1}, Lemma~\ref{lemA}, and the Cauchy-Schwarz inequality.
\medskip
\end{proof*}

\noindent
From now and until the end of this section, we shall suppose that the assumptions of Theorem~\ref{thmHHV} are fulfilled. 

\begin{lem}
\label{lemHHV1}
Let us define,
\begin{gather*}
E=\Big\{u\in H^1(\R^N); Vu^2\in L^1(\R^N) \text{ and } Au\in L^2(\R^N;\C^N)\Big\}.
\end{gather*}
Then, $H^1_{A,V}(\R^N)=E$ and $\Dr(\R^N)\inj H^1_{A,V}(\R^N)\underset{\text{dense}}{\inj}H^1(\R^N).$ In particular, each term in the integrals of $\langle\:.\:,\:.\:\rangle_{H^1_{A,V}(\R^N)}$ belongs to $L^1(\R^N).$
\end{lem}

\begin{proof*}
It is clear that $E\subset H^1_{A,V}(\R^N).$ By Lemma~\ref{lemGA}, $H^1_{A,V}(\R^N)\subset E\subset H^1(\R^N).$ It follows that $H^1_{A,V}(\R^N)=E,$ which gives the last part of the lemma, with help of Lemma~\ref{lemps}. Let $u\in H^1_{A,V}(\R^N).$ We have by \eqref{V}, \eqref{lemA1} and \eqref{thmdi1},
\begin{gather*}
\|u\|_{L^2(\R^N)}^2\le\frac1\nu\vint_{\R^N}V|u|^2\d x\le\frac1\nu\|u\|_{H^1_{A,V}(\R^N)}^2,
\end{gather*}
and
\begin{align*}
	&	\; \|\nabla u\|_{L^2(\R^N)}\le\|\nabla u+\vi Au\|_{L^2(\R^N)}+\|Au\|_{L^2(\R^N)}			\\
  \le	&	\; \|\nabla u+\vi Au\|_{L^2(\R^N)}+C \alpha_A \|\,|u|\,\|_{H^1(\R^N)}			\\
  \le	&	\; \|\nabla u+\vi Au\|_{L^2(\R^N)}+C \alpha_A\left(\|u\|_{L^2(\R^N)}+\|\nabla u+\vi Au\|_{L^2(\R^N)}\right)	\\
  \le	&	\; C\|u\|_{H^1_{A,V}(\R^N)}.
\end{align*}
Hence, $H^1_{A,V}(\R^N)\inj H^1(\R^N).$ Let $\vphi\in\Dr(\R^N).$ Let $R>0$ be such that $\supp\vphi\subset B(0,R).$ By Hölder's inequality, $V\vphi^2\in L^1(\R^N)$ and $A\vphi\in L^2(\R^N;\C^N).$ It follows that $\vphi\in H^1_{A,V}(\R^N).$ Again by Hölder's inequality and \eqref{lemGA3}, we have
\begin{gather*}
\|\vphi\|_{H^1_{A,V}(\R^N)}^2\le\|V\|_{L^1(B(0,R))}\|\vphi\|_{L^\infty(\R^N)}^2+C\|\vphi\|_{H^1(\R^N)}^2,
\end{gather*}
where $C$ does not depend on $\vphi.$ Hence, $\Dr(\R^N)\inj H^1_{A,V}(\R^N).$ Finally, since $\Dr(\R^N)\subset H^1_{A,V}(\R^N)$ and $\Dr(\R^N)\underset{\text{dense}}{\inj}H^1(\R^N),$ we conclude that $H^1_{A,V}(\R^N)\underset{\text{dense}}{\inj}H^1(\R^N).$
\medskip
\end{proof*}

\begin{lem}
\label{lemHHV2}
It holds that $\Dr(\R^N)\subset H^1_{A,V}(\R^N)$ and for any $u\in H^1_{A,V}(\R^N),$ there exists $(\vphi_n)_{n\in\N}\subset\Dr(\R^N)$ such that $\vphi_n\xrightarrow[n\to\infty]{H^1_{A,V}}u.$ In addition, we have for any $n\in\N,$ $\|\vphi_n\|_{L^p(\R^N)}\le\|u\|_{L^p(\R^N)},$ for any $p\in[1,\infty],$ as soon as $u\in L^p(\R^N).$

\end{lem}

\begin{proof*}
We  adapt the proof of Theorem~7.22, p.194, in Lieb and Loss~\cite{MR1817225}
to handle the presence of the potential $V$ in the integral $\int V|u|^2\d
x.$  By Lemma~\ref{lemGA3}, we already know that $\Dr(\R^N)\subset H^1_{A,V}(\R^N).$  Let $\xi\in C^\infty(\R;\R)$ be such that $0\le\xi\le1,$ $\xi(t)=1,$ if $|t|\le1$ and $\xi(t)=0,$ if $|t|\ge2.$ Let $n\in\N.$ Set for any $x\in\R^N,$ $\xi_n(x)=\xi\left(\frac{|x|}n\right).$ We denote by $(\rho_n)_{n\in\N}\subset\Dr(\R^N)$ any standard sequence of mollifiers.
\\
Let $u\in H^1_{A,V}(\R^N).$ Let $\eps>0.$ Let $p\in[1,\infty]$ be such that
$u\in L^p(\R^N).$  We proceed in three steps.
\\
\textbf{Step 1:} There exists $v\in H^1_{A,V}(\R^N)\cap L^\infty_\co(\R^N)$ such that $\|u-v\|_{H^1_{A,V}(\R^N)}<\frac\eps2$ and $|v|\overset{\text{a.e.}}{\le}|u|.$
\\
Let $n\in\N.$ Let for $x\in\R^N,$ $u_n(x)=\xi_n(x)\xi\left(\frac{|u(x)|}n\right)u(x).$ Then, $\supp u_n\subset\ovl B(0,2n),$ $\|u_n\|_{L^\infty(\R^N)}\le2n,$ $u_n\xrightarrow[n\to\infty]{\text{a.e. in }\R^N}u,$ $\sqrt Vu_n\xrightarrow[n\to\infty]{\text{a.e. in }\R^N}\sqrt Vu,$ $|u_n|\le|u|\in L^2(\R^N)$ and $\sqrt V|u_n|\le\sqrt V|u|\in L^2.$ It follows that $(u_n)_{n\in\N}\subset L^\infty_\co(\R^n)$ and by the dominated convergence Theorem,
\begin{gather}
\label{lemHHV21}
\lim_{n\to\infty}\vint_{\R^N}V|u-u_n|^2\d x=0 \; \text{ and } \; u_n\xrightarrow[n\to\infty]{L^2(\R^N)}u.
\end{gather}
In addition,
\begin{gather*}
\nabla u_n=\frac1n\xi^\p\left(\frac{|\:.\:|}n\right)\xi\left(\frac{|u|}n\right)\,u\,\frac{x}{|x|}+\frac1n\,\xi_n\,\xi^\p\left(\frac{|u|}n\right)u\nabla|u|
+\xi_n\xi\left(\frac{|u|}n\right)\nabla u\xrightarrow[n\to\infty]{\text{a.e. in }\R^N}\nabla u.
\end{gather*}
But, $\xi^\p\left(\frac{|u|}n\right)=0,$ if $|u|\ge2n$ so that,
\begin{gather*}
|\nabla u_n|\le\|\xi^\p\|_{L^\infty(\R)}|u|+(2\|\xi^\p\|_{L^\infty(\R)}+1)|\nabla u|\in L^2(\R^N),
\end{gather*}
by Lemmas~\ref{lemHHV1} and \ref{lemH1}. By the dominated convergence Theorem and \eqref{lemHHV21}, we then infer that $u_n\xrightarrow[n\to\infty]{H^1(\R^N)}u.$ It follows from~\eqref{lemGA3}--\eqref{lemHHV21} that $(u_n)_{n\in\N}\subset H^1_{A,V}(\R^N)$ and $u_n\xrightarrow[n\to\infty]{H^1_{A,V}(\R^N)}u.$ Pick any $n_0\in\N$ large enough to have $\|u-u_{n_0}\|_{H^1_{A,V}(\R^N)}<\frac\eps2.$ Hence the result with $v=u_{n_0}.$
\\
\textbf{Step 2:} There exists $\vphi\in\Dr(\R^N)$ such that $\|v-\vphi\|_{H^1_{A,V}(\R^N)}<\frac\eps2$ and $\|\vphi\|_{L^p(\R^N)}\le\|v\|_{L^p(\R^N)}.$
\\
Let $n\in\N.$ Let $R>1$ be such that $B(0,R)\supset\supp v.$ Set, $v_n=\rho_n\star v.$ Since $v\in H^1_c(\R^N)$ it is well-known that $v_n\in\Dr(\R^N),$ $\|v_n\|_{L^p(\R^N)}\le\|v\|_{L^p(\R^N)}$ (by Young's inequality),
\begin{align}
\label{lemHHV22}
&	\supp v_n\subset\supp\rho_n+\supp v\subset B(0,2R),	\\
\label{lemHHV23}
&	v_n\xrightarrow[n\to\infty]{H^1(\R^N)}v,
\end{align}
(see for instance Brezis~\cite{MR2759829}: Proposition~4.18, p.106; Proposition~4.20, p.107; Theorem~4.22, p.109; Lemma~9.1, p.266). Then, $v_n\in H^1_{A,V}(\R^N)$ and by \eqref{lemGA3},
\begin{gather}
\label{lemHHV24}
\nabla v_n+\vi Av_n\xrightarrow[n\to\infty]{L^2(\R^N)}\nabla v+\vi Av.
\end{gather}
By \eqref{lemHHV23}, we may extract a subsequence, that we still denote by $(v_n)_{n\in\N},$ such that $v_n\xrightarrow[n\to\infty]{\text{a.e. in }\R^N}v.$ As a consequence, $\sqrt Vv_n\xrightarrow[n\to\infty]{\text{a.e. in }\R^N}\sqrt Vv.$ Applying Young's inequality and \eqref{lemHHV22}, we see that
\begin{gather*}
\sqrt V|v_n|\le\sqrt V\|v\|_{L^\infty(\R^N)}\1_{B(0,2R)}\in L^2(\R^N).
\end{gather*}
It follows from the dominated convergence Theorem that $\vlim_{n\to\infty}\dsp\int_{\R^N}V|v-v_n|^2\d x=0,$ which gives with \eqref{lemHHV24},
\begin{gather*}
\lim_{n\to\infty}\|v-v_n\|_{H^1_{A,V}(\R^N)}=0.
\end{gather*}
We then choose $\vphi=v_{n_1},$ where $n_1\in\N$ is sufficiently large to have $\|v-v_{n_1}\|_{H^1_{A,V}(\R^N}<\frac\eps2.$ Hence Step~2.
\\
\textbf{Step 3:} Conclusion. \\
The result follows from Steps 1 and 2.
\medskip
\end{proof*}

\begin{lem}
\label{lemHHV3}
The space $H^1_{A,V}(\R^N)$ is complete.
\end{lem}

\begin{proof*}
Let $(u_n)_{n\in\N}\subset H^1_{A,V}(\R^N)$ be a Cauchy sequence. Since $H^1_{A,V}(\R^N)\inj H^1(\R^N)$ which is complete (Lemma~\ref{lemHHV1}), there exists $u\in H^1(\R^N)$ such that $u_n\xrightarrow[n\to\infty]{H^1(\R^N)}u.$  By Lemma~\ref{lemGA}, $Au\in L^2(\R^N;\C^N)$ and $\nabla u_n+\vi Au_n\xrightarrow[n\to\infty]{L^2(\R^N)}\nabla u+\vi Au.$ To conclude, it remains to show that $\sqrt Vu\in L^2(\R^N)$ and $\sqrt Vu_n\xrightarrow[n\to\infty]{L^2(\R^N)}\sqrt Vu.$ The sequence $(\sqrt Vu_n)_{n\in\N}$ being  Cauchy in $L^2(\R^N),$ it is bounded and there exists $v\in L^2(\R^N)$ such that $\sqrt Vu_n\xrightarrow[n\to\infty]{L^2(\R^N)}v.$ There exists a subsequence $(u_{n_k})_{k\in\N}\subset(u_n)_{n\in\N}$ such that $u_{n_k}\xrightarrow[k\to\infty]{\text{a.e. in }\R^N}u.$ It follows that $\sqrt Vu_{n_k}\xrightarrow[k\to\infty]{\text{a.e. in }\R^N}\sqrt Vu$ and by Fatou's Lemma, $\sqrt Vu\in L^2(\R^N).$ Let $\vphi\in\Dr(\R^N).$ By Hölder's inequality, $\sqrt V\vphi\in L^2(\R^N).$ We have for any $n\in\N,$
\begin{gather*}
\langle\sqrt Vu_n,\vphi\rangle_{L^2(\R^N),L^2(\R^N)}=\langle u_n,\sqrt V\vphi\rangle_{L^2(\R^N),L^2(\R^N)}.
\end{gather*}
By the above convergences, we can pass to the limit and we get for any $\vphi\in\Dr(\R^N),$
\begin{gather*}
\langle v,\vphi\rangle_{L^2(\R^N),L^2(\R^N)}=\langle u,\sqrt V\vphi\rangle_{L^2(\R^N),L^2(\R^N)}=\langle\sqrt V u,\vphi\rangle_{L^2(\R^N),L^2(\R^N)}.
\end{gather*}
It follows that, $v=\sqrt V u$ in $\Dr^\p(\R^N)$ and so in $L^2(\R^N).$ The lemma is proved.
\medskip
\end{proof*}

\begin{vproof}{of Theorem~\ref{thmHH}.}
By Lemma~\ref{lemHHV1}, $H^1_{A,V}(\R^N)\inj H^1(\R^N).$ It remains to show that, $H^1(\R^N)\inj H^1_{A,V}(\R^N).$ Let $u\in H^1(\R^N).$ Then, $|u|\in H^1(\R^N)$ and $\|\,|u|\,\|_{H^1(\R^N)}\le\|u\|_{H^1(\R^N)}$ (Lemma~\ref{lemH1}). By Lemma~\ref{lemA}, $Au\in L^2(\R^N;\C^N)$ and $Vu^2\in L^1(\R^N).$ As a consequence, $u\in H^1_{A,V}(\R^N)$ and by \eqref{lemV1} and \eqref{lemGA3},
\begin{gather*}
\|u\|_{H^1_{A,V}(\R^N)}\le C\|u\|_{H^1(\R^N)},
\end{gather*}
where $C$ does not depend on $u.$
\medskip
\end{vproof}

\begin{lem}
\label{lemHHV4}
The space $H^1_{A,V}(\R^N)$ is separable.
\end{lem}

\begin{proof*}
By Theorem~\ref{thmHH}, $H^1_{A,1}(\R^N)=H^1(\R^N)$ with equivalent norms. Let $H=H^1(\R^N),$ with $\|\:.\:\|_H=\|\:.\:\|_{H^1_{A,1}(\R^N)}.$ Since $(H^1(\R^N),\|\:.\:\|_{H^1(\R^N)})$ is separable, so is $(H,\|\:.\:\|_H).$ Let us define the linear operator $T$ by,
\begin{gather*}
\begin{array}{rcl}
T:H^1_{A,V}(\R^N)	&	\tends	&	L^2(\R^N)\times H	\medskip \\
			 u	& \longmapsto	&	(\sqrt Vu,u),
\end{array}
\end{gather*}
with $\|(u,v)\|_{L^2(\R^N)\times H}^2=\|u\|_{L^2(\R^N)}^2+\|v\|_H^2,$ for any $(u,v)\in L^2(\R^N)\times H.$ Clearly, $L^2(\R^N)\times H$ is separable. Thus, $T\big(H^1_{A,V}(\R^N)\big)$ is also separable (Brezis~\cite{MR2759829}: Proposition~3.25, p.73). But for any $u\in H^1_{A,V}(\R^N),$
\begin{gather*}
\|T(u)\|_{L^2(\R^N)\times H}\ge\|u\|_{H^1_{A,V}(\R^N)},
\end{gather*}
so that $H^1_{A,V}(\R^N)$ is separable.
\medskip
\end{proof*}

\begin{vproof}{of Theorem~\ref{thmHHV}.}
By Lemmas~\ref{lemHHV1}--\ref{lemHHV4}, it remains to show the continuous embeddings and the densities in \eqref{thmHHV4}. This comes from the fact that the embeddings in \eqref{thmHHV3} are dense and from the reflexivity of the spaces $\Dr(\R^N)$ and $H^1_{A,V}(\R^N).$
\medskip
\end{vproof}

\begin{vproof}{of Theorem~\ref{thmLA}.}
Estimate \eqref{thmLA1} comes from \eqref{rmkdefsol11}--\eqref{rmkdefsol13} and a straightforward calculation, while \eqref{thmLA2} is a consequence of the embeddings $H^1_{A,V}(\R^N)\inj H^1(\R^N)$ and $H^{-1}(\R^N)\inj H^{-1}_{A,V}(\R^N),$ due to \eqref{thmHHV3} and \eqref{thmHHV4}.
\medskip
\end{vproof}

\section{The set of dislocations}
\label{sod}

\begin{lem}
\label{lemAA}
Let $A$ satisfy~\eqref{A}, \eqref{Ae} and \eqref{curlA}. Then for any
$y\in\Z^N,$ there exists a unique continuous function
$\psi_y\in W^{1,N+\eps}_\loc(\R^N;\R)$ $(\psi_y\in H^1_\loc(\R;\R),$ if
$N=1)$ such that
\begin{gather}
\label{lemAApsi0}
\psi_y(0)=0,											\\
\label{psieven}
\forall x\in\R^N, \; \psi_y(x-y)+\psi_{-y}(x)=\psi_y(-y)=\psi_{-y}(y),	\\
\label{lemAApsi}
A(x+y)=A(x)+\nabla\psi_y(x),
\end{gather}
for almost every $x\in\R^N.$ In particular, $\psi_0=0$ over $\R^N.$
\end{lem}

\begin{proof*}
  Let $y\in\Z^N.$ Uniqueness for $\psi_y$ comes from \eqref{lemAApsi0} and
  \eqref{lemAApsi}, once continuity is proved. By Remark~\ref{rmkhypA2} and
  the Sobolev embedding, there exists
  $\wt{\psi_y}\in W^{1,N+\eps}_\loc(\R^N;\R)$ satisfying \eqref{lemAApsi}
  and continuous over $\R^N.$ Setting $\psi_y=\wt{\psi_y}-\wt{\psi_y}(0),$
  we see that $\psi_y$ verifies \eqref{lemAApsi0} and \eqref{lemAApsi}.
   Notice that the function $x\longmapsto0$ satisfies
  \eqref{lemAApsi} for $y=0,$ so that $\psi_0=0,$ by uniqueness. It remains
  to establish \eqref{psieven}. Applying \eqref{lemAApsi} with $y$ at the
  point $x-y$ and a second time with $-y,$ we obtain for almost every
  $x\in\R^N,$
\begin{gather*}
A(x-y)=A(x)-\nabla\psi_y(x-y)=A(x)+\nabla\psi_{-y}(x).
\end{gather*}
It follows that there exists $c\in\R$ such that,
\begin{gather*}
\forall x\in\R^N, \; \psi_y(x-y)+\psi_{-y}(x)=c.
\end{gather*}
Substituting first $x=0,$ then  $x=y$ and using \eqref{lemAApsi0} we obtain
\eqref{psieven}.
\medskip
\end{proof*}

\begin{lem}
\label{lemphi}
Let $A$ satisfy~\eqref{A}, \eqref{Ae} and \eqref{curlA}. Let
$\big(\psi_y\big)_{y\in\Z^N}$ be given by Lemma~$\ref{lemAA}.$ For any
$y\in\Z^N,$ let $\vphi_y\in W^{1,N+\eps}_\loc(\R^N;\R)$ be defined by,
\begin{gather}
\label{defphi}
\vphi_y \eqdef \psi_y-\frac12\psi_y(-y),
\end{gather}
Then $\vphi_y\in C(\R^N;\R)$ and verifies,
\begin{gather}
\label{lemphiequ}
\forall x\in\R^N, \; \vphi_y(x-y)+\vphi_{-y}(x)=0,	\\
\label{lemAphi}
A(x+y)=A(x)+\nabla\vphi_y(x),
\end{gather}
for almost every $x\in\R^N.$ Finally, $\vphi_0=0$ over $\R^N.$
\end{lem}

\begin{proof*}
By Lemma~\ref{lemAA} and \eqref{defphi}, we only have to check \eqref{lemphiequ}. The result then comes from \eqref{defphi} and \eqref{psieven}.
\medskip
\end{proof*}

\noindent
Assume that $A$ and $V$ satisfy Assumptions~\ref{assA} and \ref{assf}, respectively. For any $y\in\Z^N,$ we
define $\tau_y\in\Lr\big(H^1_{A,V}(\R^N)\big)$ as follows.
\begin{gather*}
	\begin{array}{rcl}
		\tau_y:H^1_{A,V}(\R^N)		&	\tends	&	H^1_{A,V}(\R^N)			\medskip \\
						u		&  \longmapsto	& 	e^{\vi\vphi_y}u(\:\cdot\:+y),
	\end{array}
\end{gather*}
where $\vphi_y$ is given by \eqref{defphi}. Indeed, it is clear that
$\tau_y:H^1_{A,V}(\R^N)\tends L^2(\R^N)$ is linear and,
\begin{gather}
\label{t1}
\vint_{\R^N}V|\tau_yu|^2\d x=\vint_{\R^N}V|u|^2\d x,
\end{gather}
for any $u\in H^1_{A,V}(\R^N).$ In addition, by \eqref{lemAphi}, we have for any $y\in\Z^N,$ $u\in H^1_{A,V}(\R^N)$ and almost every $x\in\R^N,$
\begin{align*}
	&	\; \nabla(\tau_yu)(x)+\vi A(x)(\tau_yu)(x)									\\
   =	&	\; \big(\nabla u(x+y)+\vi A(x)u(x+y)+\vi\nabla\vphi_y(x)u(x+y)\big)e^{\vi\vphi_y(x)}	\\
      =	&	\; \big(\nabla u(x+y)+\vi A(x+y)u(x+y)\big)e^{\vi\vphi_y(x)}.
\end{align*}
We deduce that $\tau_y:H^1_{A,V}(\R^N)\tends H^1_{A,V}(\R^N)$ is well-defined, linear and
\begin{gather*}
\|\nabla(\tau_yu)+\vi A(\tau_yu)\|_{L^2(\R^N)}=\|\nabla u+\vi Au\|_{L^2(\R^N)}.
\end{gather*}
The above estimates and \eqref{t1} permit us to see that for any $y\in\Z^N,$ $\tau_y\in\Lr\big(H^1_{A,V}(\R^N)\big)$ with $\|\tau_y\|_{\Lr(H^1_{A,V}(\R^N))}=1.$ Let
\begin{gather}
\label{D}
D \eqdef \big\{\tau_y;y\in\Z^N\big\}.
\end{gather}

\begin{prop}
\label{propDD1}
Let $A$ and $V$ satisfy Assumptions~$\ref{assA}$ and $\ref{assf},$ respectively, and let $D$ be defined by~\eqref{D}. Then $D$ is a set of unitary operators on $H^1_{A,V}(\R^N).$ In addition,
\begin{gather}
\label{propDD1-0}
\tau_0=\Id,	\\
\label{propDD1-1}
\tau_y^{-1}=\tau_y^\star=\tau_{-y},	\\
\label{propDD1-2}
\langle\tau_yu,\tau_yv \rangle_{H^1_{A,V}(\R^N)}=\langle u,v \rangle_{H^1_{A,V}(\R^N)},
\end{gather}
for any $y\in\Z^N$ and $u,v\in H^1_{A,V}(\R^N).$
\end{prop}

\begin{vproof}{of Lemma~\ref{propDD1}.}
Recall that $D$ is set of bounded linear operators on $H^1_{A,V}(\R^N).$ By Lemma~\ref{lemphi}, $\vphi_0=0$ so that $\tau_0=\Id.$ Let $y\in\Z^N$ and let $u\in H^1_{A,V}(\R^N).$ For almost every $x\in\R^N,$ one has,
\begin{gather*}
\tau_y\big(\tau_{-y}u\big)(x)=e^{\vi\vphi_y(x)}\big(\tau_{-y}u\big)(x+y)=e^{\vi\vphi_y(x)}e^{\vi\vphi_{-y}(x+y)}u(x)=u(x),
\end{gather*}
where we have used \eqref{lemphiequ} in the last equality. Still with
\eqref{lemphiequ}, we show that $\tau_{-y}\big(\tau_yu\big)=u.$ It follows that
$\tau_y$ is invertible and $\tau_y^{-1}=\tau_{-y}.$ Now, let $v\in H^1_{A,V}(\R^N).$ By a
straightforward calculation and with help of \eqref{lemphiequ} again and
\eqref{lemAphi}, we obtain
\begin{gather*}
  \langle u,\tau_y^\star v \rangle_{H^1_{A,V}(\R^N)}\eqdef \langle \tau_yu,v
  \rangle_{H^1_{A,V}(\R^N)}= \langle u,\tau_y^{-1}v \rangle_{H^1_{A,V}(\R^N)},
\end{gather*}
so that, $\tau_y^\star=\tau_y^{-1}$ which concludes the proof.
\medskip
\end{vproof}

\noindent
Let us recall the following definition (see Definition~3.1, p.60, in~Tintarev and Fieseler~\cite{MR2294665} and Proposition~\ref{propDD1}).

\begin{defi}
\label{defwD}
Let $A$ and $V$ satisfy Assumptions~$\ref{assA}$ and $\ref{assf},$ respectively, and let $D$ be defined by~\eqref{D}. Let $(u_n)_{n\in \N}\subset H^1_{A,V}(\R^N)$ and $u\in H^1_{A,V}(\R^N).$ We shall say that the sequence $(u_n)_{n\in \N}$ \textit{converges to} $u$ $D$\textit{-weakly}, which we will denote as,
\begin{gather*}
u_n\underset{n\to\infty}{\overset{D}{-\!\!\!-\!\!\!\weak}}u,
\end{gather*}
if
\begin{gather*}
\lim_{n\to\infty}\langle u_n-u,\tau_{y_n}v \rangle_{H^1_{A,V}(\R^N)}=0,
\end{gather*}
for any sequence $(\tau_{y_n})_{n\in\N}\subset D$ and $v\in H^1_{A,V}(\R^N).$
\end{defi}

\begin{nota}
\label{notawD}
Let $A$ and $V$ satisfy Assumptions~$\ref{assA}$ and $\ref{assf},$ respectively, and let $D$ be defined by~\eqref{D}. Let $(\tau_{y_n})_{n\in\N}\subset D.$ We shall write,
\begin{gather*}
\tau_{y_n}\underset{n\to\infty}{-\!\!\!-\!\!\!\weak}0,
\end{gather*}
to mean that for any $u\in H^1_{A,V}(\R^N),$ $\tau_{y_n}u\underset{n\to\infty}{-\!\!\!-\!\!\!\weak}0$ in $H^1_{A,V}$-weakly, or equivalently,
\begin{gather*}
\lim_{n\to\infty}\langle\tau_{y_n}u_,v \rangle_{H^1_{A,V}(\R^N)}=0,
\end{gather*}
for any $u,v\in H^1_{A,V}(\R^N).$
\end{nota}

\begin{rmk}
\label{rmkdefnotawD}
If $u_n\underset{n\to\infty}{\overset{D}{-\!\!\!-\!\!\!\weak}}0$ then $u_n\underset{n\to\infty}{-\!\!\!-\!\!\!\weak}0$ in $H^1_{A,V}$-weakly. In particular, $u_n\underset{n\to\infty}{\overset{H^1_\w}{-\!\!\!-\!\!\!\weak}}0$ and for any sequence $(\tau_{y_n})_{n\in\N}\subset D,$ $\tau_{y_n}u_n\underset{n\to\infty}{\overset{H^1_\w}{-\!\!\!-\!\!\!\weak}}0.$ Indeed, this follows from \eqref{propDD1-0}, \eqref{propDD1-1} and \eqref{thmHHV3}.
\end{rmk}

\begin{lem}
\label{lemD}
Let $A$ and $V$ satisfy Assumptions~$\ref{assA}$ and $\ref{assf},$ respectively. Let $(y_k)_k,(z_k)_k\subset\Z^N.$ Then,
\begin{gather}
\label{lemD1}
\tau_{y_k}\tau_{z_k}\underset{k\to\infty}{-\!\!\!-\!\!\!\weak}0 \iff |y_k+z_k|\xrightarrow[k\to\infty]{}\infty.
\end{gather}
Moreover, if $\tau_{y_k}\tau_{z_k}\cancel{\underset{k\to\infty}{-\!\!\!-\!\!\!\weak}}0$ then $\big(\tau_{y_k}\tau_{z_k}\big)_k$ admits a constant subsequence.
\end{lem}

\begin{proof*}
Let $(y_k)_k,(z_k)_k\subset\Z^N.$
\\
\textbf{Step~1:} If $\vliminf_{k\to\infty}|y_k+z_k|<\infty$ then $(y_k+z_k)_k$ admits a constant subsequence.
\\
Indeed, if $\vliminf_{k\to\infty}|y_k+z_k|<\infty$ then $(y_k+z_k)_k$ admits a bounded subsequence, from which we extract a convergent subsequence $\big(y_{k_\ell}+z_{k_\ell}\big)_\ell.$ Since $\big(y_{k_\ell}+z_{k_\ell}\big)_\ell$ converges in $\Z^N,$ Step~1 follows.
\\
\textbf{Step~2:} Proof of $\Longrightarrow.$
\\
We show the contraposition. Assume that $\vliminf_{k\to\infty}|y_k+z_k|<\infty.$ By Step~1, there exists $\big(y_{k_\ell}+z_{k_\ell}\big)_\ell\subset(y_k+z_k)_k$ such that for any $\ell\in\N,$ $y_{k_\ell}+z_{k_\ell}=y_{k_1}+z_{k_1}.$ Let $u\in H^1_{A,V}(\R^N)\setminus\{0\}$ and $v=\tau_{y_{k_1}}\tau_{z_{k_1}}u.$ It follows from \eqref{propDD1-2} that,
\begin{gather*}
\forall\ell\in\N, \; \langle\tau_{y_{k_\ell}}\tau_{z_{k_\ell}}u,v\rangle_{H^1_{A,V}(\R^N)}=\|u\|_{H^1_{A,V}(\R^N)}^2>0,
\end{gather*}
and so, $\tau_{y_k}\cancel{\underset{k\to\infty}{-\!\!\!-\!\!\!\weak}}0.$
\\
\textbf{Step~3:} Proof of $\Longleftarrow.$
\\
Assume $|y_k+z_k|\xrightarrow[k\to\infty]{}\infty.$ Let $\vphi,\psi\in\Dr(\R^N).$ Then for any $k\in\N$ large enough, $\supp(\tau_{y_k}\tau_{z_k}\vphi)\cap\supp\psi=\emptyset,$ so that,
\begin{gather}
\label{prooflemD}
\langle\tau_{y_k}\tau_{z_k}\vphi,\psi \rangle_{H^1_{A,V}(\R^N)}\xrightarrow{k\to\infty}0.
\end{gather}
Let $u,v\in H^1_{A,V}(\R^N).$ Let $\eps>0.$ By density (Theorem~\ref{thmHHV}), there exists $(\vphi_n)_n,(\psi_n)_n\subset\Dr(\R^N)$ such that, $\vphi_n\xrightarrow[n\to\infty]{H^1_{A,V}(\R^N)}u$ and $\psi_n\xrightarrow[n\to\infty]{H^1_{A,V}(\R^N)}v.$ Let $n_0\in\N$ be such that,
\begin{gather*}
\|v\|_{H^1_{A,V}(\R^N)}\|u-\vphi_{n_0}\|_{H^1_{A,V}(\R^N)}+\|\vphi_{n_0}\|_{H^1_{A,V}(\R^N)}\|v-\psi_{n_0}\|_{H^1_{A,V}(\R^N)}\le\eps,
\end{gather*}
for any $n\ge n_0.$ We then infer with help of~\eqref{propDD1-2} and Cauchy-Schwarz's inequality that for any $k\in\N,$
\begin{align*}
	&	\; |\langle\tau_{y_k}\tau_{z_k}u,v \rangle_{H^1_{A,V}}|\le| \langle\tau_{y_k}\tau_{z_k}(u-\vphi_{n_0}),v \rangle_{H^1_{A,V}}|+|\langle
          		\tau_{y_k}\tau_{z_k}\vphi_{n_0},v-\psi_{n_0} \rangle_{H^1_{A,V}}|+| \langle\tau_{y_k}\tau_{z_k}\vphi_{n_0},\psi_{n_0} \rangle_{H^1_{A,V}}|\\
  \le	&	\; \|v\|_{H^1_A}\|u-\vphi_{n_0}\|_{H^1_{A,V}}+\|\vphi_{n_0}\|_{H^1_A}\|v-\psi_{n_0}\|_{H^1_{A,V}}
			+|\langle\tau_{y_k}\tau_{z_k}\vphi_{n_0},\psi_{n_0} \rangle_{H^1_{A,V}}|													\\
  \le	&	\; \eps+| \langle\tau_{y_k}\tau_{z_k}\vphi_{n_0},\psi_{n_0} \rangle_{H^1_{A,V}}|.
\end{align*}
By~\eqref{prooflemD}, if follows that: $\vlimsup_{k\to\infty}|\langle\tau_{y_k}\tau_{z_k}u,v\rangle_{H^1_{A,V}(\R^N)}|\le\eps.$ Since $\eps>0$ is arbitrary,
we then get that for any $u,v\in H^1_{A,V}(\R^N),$ $\big\langle\tau_{y_k}\tau_{z_k}u,v\big\rangle_{H^1_{A,V}(\R^N)}\xrightarrow{k\to\infty}0,$ which is the desired result.
\\
\textbf{Step~4:} If $\tau_{y_k}\tau_{z_k}\cancel{\underset{k\to\infty}{-\!\!\!-\!\!\!\weak}}0$ then $\big(\tau_{y_k}\tau_{z_k}\big)_k$ admits a constant subsequence.
\\
Now assume that $\tau_{y_k}\tau_{z_k}\cancel{\underset{k\to\infty}{-\!\!\!-\!\!\!\weak}}0.$ By \eqref{lemD1}, this means $\vliminf_{k\to\infty}|y_k+z_k|<\infty,$ and we conclude with help of Step~1.
\medskip
\end{proof*}

\begin{prop}
\label{propDD2}
Let $A$ and $V$ satisfy Assumptions~$\ref{assA}$ and $\ref{assf},$ respectively, and let $D$ be defined by~\eqref{D}. Then $D$ is a set of dislocations on $H^1_{A,V}(\R^N).$
\end{prop}

\begin{proof*}
By Proposition~3.1 p.61 in Tintarev and Fieseler~\cite{MR2294665}, it is
sufficient to show that if $(y_k)_k\subset\Z^N$ is such that
$\tau_{y_k}\cancel{\underset{k\to\infty}{-\!\!\!-\!\!\!\weak}}0$ then $\tau_{y_k}$ has a strongly
convergence subsequence. This is a consequence of~\eqref{propDD1-0} and  Lemma~\ref{lemD}.
\medskip
\end{proof*}

\section{Cocompactness}
\label{coc}

\begin{thm}
\label{thmcoc}
Let $A$ and $V$ satisfy Assumptions~$\ref{assA}$ and $\ref{assf},$ respectively, and let $D$ be defined by~\eqref{D}. Let $(u_k)_{k\in\N}$ be a bounded sequence in $H^1_{A,V}(\R^N).$ Let $p\in(2,\2)$ $(p\in(2,\infty],$ if $N=1).$ Then $(u_k)_{k\in\N}$ is bounded in $H^1(\R^N)$ and we have the
following result.
\begin{gather*}
u_k\underset{k\to\infty}{\overset{D}{-\!\!\!-\!\!\!\weak}}0 \iff u_k\xrightarrow[k\to\infty]{L^p(\R^N)}0.
\end{gather*}
\end{thm}

\begin{proof*}
Let $(u_k)_{k\in\N}$ be a bounded sequence in $H^1_{A,V}(\R^N)$ and let $p$ be as in the theorem. By \eqref{thmHHV3} and Sobolev' embedding, $(u_k)_{k\in\N}$ is bounded in $H^1(\R^N)$ and so in $L^p(\R^N).$ Assume that $u_k\underset{k\to\infty}{\overset{D}{-\!\!\!-\!\!\!\weak}}0.$ By Remark~\ref{rmkdefnotawD}, $\tau_{y_k}u_k\underset{k\to\infty}{\overset{H^1_\w}{-\!\!\!-\!\!\!\weak}}0,$ for any $(\tau_{y_k})_k\subset D.$ Suppose $p<\infty.$ We claim that,
\begin{gather}
\label{demthmcoc}
\forall k\in\N, \; \exists y_k\in\Z^N \text{ such that } \sup_{y\in\Z^N}\vint_{Q-y}|u_k|^p\d x=\vint_Q|\tau_{y_k}u_k|^p\d x.
\end{gather}
Indeed, if $\vsup_{y\in\Z^N}\vint_{Q-y}|u_k|^p\d x=0,$ there is nothing to
prove. If $\vsup_{y\in\Z^N}\vint_{Q-y}|u_k|^p\d x= \delta >0$ then if the
supremum in $y$ was not a maximum then there would be an infinite number of
$y \in \Z^N$ such that $\vint_{Q-y}|u_k|^p\d x > \frac\delta2,$
contradicting the fact that $(u_k)_k$ is bounded in $L^p(\R^N).$
 \\
By the Sobolev embedding $H^1(Q)\inj L^p(Q)$ and translation, there exists
$C>0$ such that for any $k\in\N$ and $y\in\Z^N,$ $\|u_k\|_{L^p(Q-y)}^2\le
C\|u_k\|_{H^1(Q-y)}^2.$ Multiplying the both sides by
$\|u\|_{L^p(Q-y)}^{p-2},$ we get
\begin{gather*}
\vint_{Q-y}|u_k|^p\d x\le C\|u_k\|_{H^1(Q-y)}^2\left(\;\vint_{Q-y}|u_k|^p\d x\right)^{\frac{p-2}p}.
\end{gather*}
Summing over $y\in\Z^N,$ we obtain for any $k\in\N,$
\begin{gather*}
\|u_k\|_{L^p(\R^N)}^p\le C\|u_k\|_{H^1(\R^N)}^2\sup_{y\in\Z^N}\left(\;\vint_{Q-y}|u_k|^p\d x\right)^{\frac{p-2}p}.
\end{gather*}
For any $k\in\N,$ let $y_k\in\Z^N$ be given by~\eqref{demthmcoc}. Noticing
that $\vsup_{k\in\N}\|u_k\|_{H^1(\R^N)}<\infty,$ we infer from the
compactness of the Sobolev embedding $H^1(Q)\inj L^p(Q)$ that
\begin{gather*}
\forall k\in\N, \; \|u_k\|_{L^p(\R^N)}^p\le C\|\tau_{y_k}u_k\|_{L^p(Q)}^{p-2}\xrightarrow{k\to\infty}0,
\end{gather*}
since $\tau_{y_k}u_k\underset{k\to\infty}{\overset{H^1_\w}{-\!\!\!-\!\!\!\weak}}0.$ When $N=1$ and $p=\infty,$ we
use the above result and Gagliardo-Nirenberg's inequality to see that,
\begin{gather*}
\|u_k\|_{L^\infty(\R)}\le C\|u_k\|_{L^4(\R)}^\frac23\|u_k\|_{H^1(\R)}^\frac13\le C\|u_k\|_{L^4(\R)}^\frac23\xrightarrow{k\to\infty}0.
\end{gather*}
Let us prove the converse and assume that $u_k\xrightarrow[k\to\infty]{L^p(\R^N)}0.$ Note that if $N=1$ and $p=\infty$ then,
\begin{gather*}
\|u_k\|_{L^4(\R)}^2\le\|u_k\|_{L^2(\R)}\|u_k\|_{L^\infty(\R)}\le C\|u_k\|_{L^\infty(\R)}\xrightarrow{k\to\infty}0.
\end{gather*}
So we may assume that $p<\infty.$ Let $(\tau_{y_k})_k\in D.$ Since for any $k\in\N,$ $\|\tau_{y_k}u_k\|_{L^p(\R^N)}=\|u_k\|_{L^p(\R^N)}$ and $\|\tau_{y_k}u_k\|_{H^1_{A,V}(\R^N)}=\|u_k\|_{H^1_{A,V}(\R^N)}$ by~\eqref{propDD1-2}, we obtain that for some $\big(\tau_{y_{k_\ell}}\big)_\ell\subset(\tau_{y_k})_k$ and $u\in H^1_{A,V}(\R^N),$
\begin{gather*}
\tau_{y_k}u_k\xrightarrow[k\to\infty]{}0, \text{ in } L^p(\R^N),								\\
\tau_{y_{k_\ell}}u_{k_\ell}\underset{\ell\to\infty}{-\!\!\!\weak}u, \text{ in } H^1_{A,V}\text{-weakly.}
\end{gather*}
In particular, both convergences hold in $\Dr^\p(\R^N),$ so that $u=0$ and $\tau_{y_k}u_k\underset{k\to\infty}{-\!\!\!\weak}0,$ in $H^1_{A,V}$-weakly, for the whole sequence $(\tau_{y_k}u_k)_k.$ This concludes the proof.
\medskip
\end{proof*}

\section{An associated critical value function and proof of the main result}
\label{scvf}

\begin{prop}
\label{propg}
Let $N\ge1.$ Let $g$ and $F$ be as in Assumption~$\ref{assf},$ where $f$ satisfies~\eqref{f1} and \eqref{f2} and let $\psi$ be defined by~\eqref{psi}. Then the following holds.
\begin{enumerate}
\item
\label{propg1}
$\psi\in C^1(H^1(\R^N);\R),$ $\psi^\p=g$ and $\psi$ and $\psi^\p$ are bounded on bounded sets.
\item
\label{propg2}
$\forall(u,v)\in H^1(\R^N)\times H^1(\R^N),$ $\langle g(u),v\rangle_{H^{-1}(\R^N),H^1(\R^N)}=\Re\dsp\vint_{\R^N}g(u)(x)\ovl{v(x)}\d x.$
\item
\label{propg3}
Let $(u_n)_n,(v_n)_n \subset H^1(\R^N)$ be bounded. If $\vlim_{n\to\infty}\|u_n-v_n\|_{L^p(\R^N)}=0,$ for some $p\in[1,\infty],$ then $\vlim_{n\to\infty}|\psi(u_n)-\psi(v_n)|=0.$
\item
\label{propg4}
Let $u\in H^1(\R^N).$ If $u_n\underset{n\to\infty}{\overset{H^1_\w}{-\!\!\!-\!\!\!-\!\!\!\weak}}u$ then $g(u_n)\underset{n\to\infty}{\overset{H^{-1}_\w}{-\!\!\!-\!\!\!-\!\!\!\weak}}g(u).$
\end{enumerate}
\end{prop}

\noindent
Proposition~\ref{propg} is well-known but with some slightly different assumptions on $f$ and, in all cases, for real-valued functions. It can be adapted and for the convenience of the reader, we postpone its proof to the Appendix~\ref{appendixA}.
 
\begin{rmk}
\label{rmkpropg}
If $N\ge3$ then, under the hypotheses of Proposition~\ref{propg}, the conclusions may be slightly more general as follows. We first recall that,
\begin{gather*}
\Dr(\R^N)\inj E\eqdef L^2(\R^N)\cap L^\2(\R^N) \; \text{ with dense embedding},		\\
E^\star=L^2(\R^N)+L^{\2^\p}(\R^N)\inj\Dr^\p(\R^N) \; \text{ with dense embedding},
\end{gather*}
where, $\|u\|_E=\|u\|_{L^2(\R^N)}+\|u\|_{L^\2(\R^N)}.$ See, for instance, Bergh and Löfström~\cite{MR0482275} (Lemma~2.3.1, p.24--25 and Theorem~2.7.1, p.32). Then the following holds.
\begin{enumerate}
\item
\label{rmkpropg1}
$\psi\in C^1(E;\R),$ $\psi^\p=g\in C(E;E^\star)$ and $\psi$ and $\psi^\p$ are bounded on bounded sets.
\item
\label{rmkpropg2}
$\forall(u,v)\in E\times E,$ $\langle g(u),v\rangle_{E^\star,E}=\Re\dsp\vint_{\R^N}g(u)(x)\ovl{v(x)}\d x.$
\item
\label{rmkpropg3}
Let $(u_n)_n,(v_n)_n \subset E$ be bounded. If $\vlim_{n\to\infty}\|u_n-v_n\|_{L^p(\R^N)}=0,$ for some $p\in[1,\infty],$ then $\vlim_{n\to\infty}|\psi(u_n)-\psi(v_n)|=0.$
\item
\label{rmkpropg4}
Let $u\in H^1(\R^N).$ If $u_n\underset{n\to\infty}{\overset{H^1_\w}{-\!\!\!-\!\!\!-\!\!\!\weak}}u$ then $g(u_n)\underset{n\to\infty}{\overset{E^\star_\w}{-\!\!\!-\!\!\!-\!\!\!\weak}}g(u).$
\end{enumerate}
For more details, see the proof of Proposition~\ref{propg} in the Appendix~\ref{appendixA}.
\end{rmk}

\noindent
From now and until the end of this section, we shall suppose that Assumptions~$\ref{assA}$ and $\ref{assf}$ are fulfilled. In particular, by~\eqref{psi+}, $g\not\equiv0.$ Note that when $g\equiv0$ then by Remark~\ref{rmkthmmain}, $u\equiv0$ is the unique solution to \eqref{nlsp}.
\medskip \\
Let for any $t\ge0,$ $\vsS_t \eqdef \Big\{u \in\Dr(\R^N); \|u\|_{H^1_{A,V}(\R^N)}^2 = t \Big\},$ $\ovl{\vsS_t}\eqdef \Big\{u \in H^1_{A,V}(\R^N); \|u\|_{H^1_{A,V}(\R^N)}^2 = t \Big\},$
$\ovl{B_t} \eqdef \Big\{u \in H^1_{A,V}(\R^N); \|u\|_{H^1_{A,V}(\R^N)}^2 \le t\Big\}$ and
\begin{gather}
\label{eqcvf}
\gamma(t) \eqdef \sup_{u \in\ovl{\vsS_t}} \psi(u)=\sup_{u \in \vsS_t} \psi(u),
\end{gather}
where the second equality in~\eqref{eqcvf} comes from density of $\Dr(\R^N)$ in $H^1_{A,V}(\R^N)$ (Theorem~\ref{thmHHV}). Furthermore let
\begin{gather}
\label{interval}
I_\gamma \eqdef \left(2\inf_{t \neq s} \frac{\gamma(t)-\gamma(s)}{t-s},2\sup_{t \neq s} \frac{\gamma(t)-\gamma(s)}{t-s} \right),
\end{gather}
and for any $\rho>0,$
\begin{gather}
\label{Grho}
\forall u\in H^1_{A,V}(\R^N), \; G_\rho (u) \eqdef \frac{\rho}{2}\|u\|_{H^1_{A,V}(\R^N)}^2 - \psi(u),	\\
\label{Gammarho}
\forall t\ge0, \; \Gamma_\rho (t) \eqdef \frac{\rho}{2} t-\gamma(t).
\end{gather}
Note that by Proposition~\ref{propg} and Theorems~\ref{thmHHV} and \ref{thmLA}, $G_\rho\in C^1(H^1_{A,V}(\R^N);\R)$ and
\begin{gather*}
G^\p_\rho(u)=\rho(-\Delta_A u+Vu)-g(u), \; \text{ in } \; H^{-1}_{A,V}(\R^N).
\end{gather*}
It follows that for $\lambda = \dfrac1\rho,$ $u_\rho$ is a weak solution to \eqref{nlsp} if, and only if, $G^\p_\rho(u _\rho) = 0.$

\begin{lem}
\label{lemLions}
The function $\gamma$ defined by~\eqref{eqcvf} is continuous and nondecreasing over $[0,\infty)$ and is locally Lipschitz continuous over $(0,\infty).$ Furthermore, $\gamma$ has a derivative at $t=0$ and $\gamma^\p(0)=0.$ In addition, for any $a\ge0$ and $b\ge0,$
\begin{equation}
\label{lions}
\gamma (a) + \gamma (b) \le \gamma(a+b). 
\end{equation}
Finally, $I_\gamma\neq\emptyset$ and $I_\gamma=\left(0,2 \vsup_{t \neq s} \frac{\gamma(t)-\gamma(s)}{t-s} \right).$
\end{lem}

\begin{proof*}
Let $u \in H^1_{A,V}(\R^N).$
\\
Let $\theta >0.$ Let $(v_k)_{k\in\N} \subset \vsS_1$ be such that $v_k \overset{}{-\!\!\!-\!\!\!-\!\!\!\weak} 0$ in $H^1_{A,V}$-weakly and $\supp v_k \subset Q.$ By Theorem~\ref{thmHHV} and compactness, $v_k \underset{k\to\infty}{\overset{H^1_\w}{-\!\!\!-\!\!\!-\!\!\!\weak}} 0$ and $v_k\xrightarrow[k\to\infty]{L^p(\R^N)}0,$ for any $p\in(2,\2),$ It follows from Property~\ref{propg3} of Proposition~\ref{propg} that,
\begin{gather}
\label{prooflem510}
\psi(u + \theta v_k)\xrightarrow{k\to\infty}\psi(u) \; \text{ and } \;
\|u+\theta v_k\|_{H^1_{A,V}(\R^N)}^2\xrightarrow{k\to\infty}\|u\|_{H^1_{A,V}(\R^N)}^2+\theta^2.
\end{gather}
Let $t>0.$ Let $(u_k)_{k\in\N} \subset \vsS_t$ be a such that $\psi(u_k)\xrightarrow[]{k\to\infty}\gamma(t)$ and $\supp u_k\subset  B(0,R_k).$ Since $\Dr(\R^N)$ is dense in $H^1_{A,V}(\R^N)$ (Theorem~\ref{thmHHV}), we may find $(w_k)_{k\in\N} \subset \Dr(\R^N)$ such that $\supp w_k\subset  B(0,r_k)$ and $w_k\xrightarrow[k\to\infty]{H^1_{A,V}(\R^N)}u.$ Let  $(y_k)_k \subset \Z^N$ with $|y_k|>R_k+r_k$ and let $v_k=\tau_{y_k}u_k\in H^1_{A,V}(\R^N).$ It follows that,
\begin{gather}
\label{prooflem511}
\forall k\in\N, \; \supp v_k\cap\supp w_k=\emptyset,
\end{gather}
from which we deduce for any $k\in\N,$ $\psi(v_k+w_k)=\psi(v_k)+\psi(w_k).$ By Theorem~\ref{thmHHV}, Proposition~\ref{propg} and the fact that $\psi$ is invariant with respect to $D,$ we have for any $k\in\N,$
\begin{align*}
	&	\; \big|\psi(u+v_k)-\big(\psi(u)+\gamma(t)\big)\big|									\\
  \le	&	\; \big|\psi(u+v_k)-\psi(w_k+v_k)\big|+|\psi(w_k)-\psi(u)|+|\psi(v_k)-\gamma(t)|				\\
  \le	&	\; C\|w_k-u\|_{H^1_{A,V}(\R^N)}+|\psi(u_k)-\gamma(t)|\xrightarrow{k\to\infty}0.
\end{align*}
It follows that,
\begin{gather}
\label{prooflem512}
\lim_{k\to\infty}\psi(u+v_k)=\psi(u)+\gamma(t).
\end{gather}
Finally, by \eqref{prooflem511} and Cauchy-Schwarz's inequality,
\begin{gather*}
|\langle u,v_k\rangle_{H^1_{A,V}(\R^N)}=|\langle u-w_k,v_k\rangle_{H^1_{A,V}(\R^N)}|\le\sqrt t\|w_k-u\|_{H^1_{A,V}(\R^N)}\xrightarrow{k\to\infty}0,
\end{gather*}
from which we get with help of \eqref{propDD1-2},
\begin{gather}
\label{prooflem513}
\lim_{k\to\infty}\|u+v_k\|_{H^1_{A,V}(\R^N)}^2=\|u\|_{H^1_{A,V}(\R^N)}^2+t.
\end{gather}
And since $\psi^\p$ is bounded on bounded sets (Proposition~\ref{propg}), we conclude with \eqref{thmHHV3} and \eqref{thmHHV4} that there exists $C_t>0$ such that,
\begin{gather}
\label{prooflem514}
\forall u \in \ovl{B_t}, \; \left|\langle \psi^\p (u), u \rangle_{H^{-1}_{A,V}(\R^N),H^1_{A,V}(\R^N)}\right|\le C_t.
\end{gather}
By \eqref{prooflem510}, \eqref{prooflem512}--\eqref{prooflem514} and \cite[Theorem 2.1]{MR1306593}, it follows that $\gamma$ is locally Lipschitz continuous and nondecreasing over $(0,\infty)$ and \eqref{lions} holds true. Now, let us prove that $\gamma^\p(0)=\vinf_{t \neq s}\frac{\gamma(t)-\gamma(s)}{t-s}=0.$ Let $\eps >0.$ We let $\kappa=1,$ if $N\ge3$ and $\kappa=0,$ if $N=2.$ Let $C_\eps>0$ and $p_\eps>2$ be given by \eqref{f1}--\eqref{f2}. By \eqref{thmHHV3} and the Sobolev embeddings, there exist $C>0,$ which does not depend on $\eps,$ and $C^\p_\eps>0$ such that for any $t>0,$
\begin{align*}
0\le	&	\; \inf_{t \neq s}\frac{\gamma(t)-\gamma(s)}{t-s}\le\frac{\gamma(t)-\gamma(0)}t\le\frac1t\sup_{u \in \vsS_t}\int_{\R^N}|F(x,|u|)|\d x	\\
  \le	&	\; \frac1t\sup_{u \in \vsS_t}\left(\eps \int_{\R^N}(|u|^2 +\kappa|u|^\2)\d x+C_\eps\int_{\R^N}|u|^{p_\eps}\d x\right)					\\
  \le	&	\; \frac{C\eps}t\left(\sup_{u \in \vsS_t}\|u\|_{H^1_{A,V}(\R^N)}^2+\kappa\sup_{u \in \vsS_t}\|u\|_{H^1_{A,V}(\R^N)}^\2\right)
			+\frac{C^\p_\eps C_\eps}t\sup_{u \in \vsS_t}\|u\|_{H^1_{A,V}(\R^N)}^{p_\eps}										\\
  \le	&	\; C\eps\left(1+\kappa t^\frac{\2-2}2\right)+C^\p_\eps C_\eps t^\frac{p_\eps-2}2\xrightarrow{t\searrow0}C\eps.
\end{align*}
Since $\eps$ is arbitrary, we can conclude that $\vinf_{t \neq s}\frac{\gamma(t)-\gamma(s)}{t-s}=\gamma^\p(0)=\vlim_{t\searrow0}\frac{\gamma(t)-\gamma(0)}t=0.$ Finally, if $I_\gamma$ where empty then we would have for any $t\ge0,$ $\gamma(t)=0.$ But this would yield $\psi(u)\le0,$ for any $u\in H^1(\R^N),$ contradicting \eqref{psi+}.
\medskip
\end{proof*}

\noindent
We shall use the well-known following result.

\begin{thm}[\textbf{\cite{MR2294665}, Theorem 3.1, p.62-63}]
\label{thmTin}
Let $(u_k)_{k\in\N}\subset H^1_{A,V}(\R^N)$ be a bounded sequence and let $t_0\eqdef\vlimsup_{k\to\infty}\|u_k\|_{H^1_{A,V}(\R^N)}^2.$ Then, up to subsequence that we will still denote by $(u_k)_{k\in\N},$ there exist $\D\subset \N,$ $(w^n)_{n\in\D} \subset H^1_{A,V}(\R^N)$ and $\left(\tau_{y_k^n}\right)_{(k,n)\in\N\times\D} \subset D$ such that for any $(n,m)\in\D\times\D,$
\begin{align}
\label{weak_limit}
&	\;\; \tau_{-y_k^n}u_k\underset{k\to\infty}{-\!\!\!-\!\!\!-\!\!\!\weak} w^n, \; \text{ in } H^1_{A,V}\text{-weakly,}		\\
\label{asymp_ortho}
&	\;\; \lim_{k\to\infty}|y_k^m-y_k^n|=\infty, \; \text{ for } n \neq m,									\\
\label{norm_bound}
&	\;\; \sum_{n \in \D} \|w^n\|_{H^1_{A,V}(\R^N)}^2 \le t_0,										\\
\label{D_weak_limit}
&	\;\; u_k - \sum_{n \in \D} \tau_{y_k^n} w^n \underset{k\to\infty}{\overset{D}{-\!\!\!-\!\!\!-\!\!\!\weak}}0,
\end{align}
where the series in~\eqref{D_weak_limit} converges uniformly in $k\in\N.$ 
\end{thm}

\begin{proof*}
Since $H^1_{A,V}(\R^N)$ is a separable infinite-dimensional Hilbert space and $D$ is a set of dislocations on $H^1_{A,V}(\R^N)$ (\eqref{thmHHV2} and Proposition~\ref{propDD2}) and using \eqref{propDD1-1}, we may apply \cite[Theorem 3.1, p.62-63]{MR2294665} which asserts, up to subsequence that we will still denote by $(u_k)_{k\in\N},$ the existence of $\D\subset \N,$ $(w^n)_{n\in\D} \subset H^1_{A,V}(\R^N)$ and $\left(\tau_{y_k^n}\right)_{(k,n)\in\N\times\D} \subset\D$ satisfying \eqref{weak_limit}, \eqref{norm_bound}, \eqref{D_weak_limit} and $\tau_{-y_k^n}\tau_{y_k^m} \underset{k\to\infty}{-\!\!\!-\!\!\!-\!\!\!\weak}0,$ for  $n \neq m.$ This last estimate and Lemma~\ref{lemD} yields \eqref{asymp_ortho}.
\medskip
\end{proof*}

\begin{lem}
\label{lemPS}
For almost every \(\rho \in I_\gamma,\) there exist $c(\rho)>0$ and a bounded critical sequence \((u_k)_{k\in\N}\subset H^1_{A,V}(\R^N)\) that is,
\begin{gather}
\label{eqPSB}
(u_k)_{k\in\N}\subset H^1_{A,V}(\R^N) \text{ is bounded,}			\\
\label{eqPS}
\begin{cases}
G_\rho (u_k) \xrightarrow{k\to\infty} c(\rho)>0,					\medskip \\
G_\rho^\p(u_k) \xrightarrow[k\to\infty]{H^{-1}_{A,V}(\R^N)}0.
\end{cases}
\end{gather}
In addition, for every $\rho\in I_\gamma,$ there exist $c(\rho)>0$ and a sequence \((u_k)_{k\in\N}\subset H^1_{A,V}(\R^N)\) satisfying \eqref{eqPS}.
\end{lem}

\begin{proof*}
The proof of \cite[Theorem 2.15]{MR1306593} can be adapted to prove Lemma \ref{lemPS}. Let $\rho_0 \in  I_\gamma.$ Then \(\Gamma_{\rho_0}\) is not monotone nondecreasing. Indeed, if so then for any \(t_1 < t_2\) we would have
\begin{gather*}
\frac{\rho_0}2 t_1 - \gamma(t_1) \le \frac{\rho_0}2 t_2 -\gamma (t_2),
\end{gather*}
which  implies $2\vsup_{t \neq s} \frac{\gamma(t)-\gamma(s)}{t-s}\le\rho_0,$
contradicting the fact that $\rho_0\in I_\gamma.$ (Similarly \(\Gamma_{\rho_0}\) is not monotone nonincreasing.) Therefore, we can find $0<t_0<t_1$ and a $\delta>0$ such that $\Gamma_{\rho_0}(t_0)>\Gamma_{\rho_0}(t_1)+3\delta>3\delta$ (we recall that by Lemma~\ref{lemLions}, $\Gamma_{\rho_0}(0)=0$ and $\Gamma_{\rho_0}^\p(0)>0).$ Also, it is clear that  the mapping $\rho\longmapsto\Gamma_\rho(t_0)-\Gamma_\rho(t_1)$ is continuous over $[0,\infty)$ so that there exists $\delta_0(\rho_0)>0$ such that for any $\rho\in I_{\rho_0}\eqdef\big(\rho_0-\delta_0(\rho_0),\rho_0+\delta_0(\rho_0)\big),$ $\Gamma_\rho(t_0)>\Gamma_\rho(t_1)+2\delta>2\delta.$ But it follows from the definition of \(\gamma\) that there is a \(u_1 \in \vsS_{t_1}\) such that $\psi(u_1)>\gamma(t_1)-\delta.$ Thus, for any $\rho\in I_{\rho_0}$ and any \(u \in \ovl{\vsS_{t_0}},\)
\begin{gather}
\label{prooflemmax}
G_\rho (u) \ge \Gamma_\rho (t_0)>\Gamma_\rho(t_1)+2\delta>G_\rho(u_1)+\delta.
\end{gather}
Denoting by $\Lambda\eqdef\left\{\xi\in C\left([0,1];H^1_{A,V}(\R^N)\right); \; \xi(0)=0 \text{ and } \xi(1)=u_1\right\},$ it follows from \eqref{prooflemmax} that the following holds.
\begin{gather}
\label{MPG}
\begin{cases}
\text{For any } \rho_0\in I_\gamma, \text{ there exist } \delta_0(\rho_0)>0 \text{ and } u_1\in H^1_{A,V}(\R^N)\setminus\{0\}		\medskip \\
\text{such that for any } \rho\in I_{\rho_0}\eqdef\big(\rho_0-\delta_0(\rho_0),\rho_0+\delta_0(\rho_0)\big),					\medskip \\
c(\rho)\eqdef\dsp\inf_{\xi\in\Lambda}\max_{t\in[0,1]}G_\rho\big(\xi(t)\big)>G_\rho(u_1)>G_\rho(0).
\end{cases}
\end{gather}
Thus \(G_{\rho_0}\) has mountain pass geometry and we can find a critical sequence satisfying \eqref{eqPS} by the Mountain Pass Theorem (see, for instance, \cite[Theorem 6.2, p.144]{MR2294665}). Now, let us show that for almost every $\rho\in I_\gamma,$ there exists a bounded critical sequence. As we shall see, this is almost a direct consequence of \eqref{MPG} and \cite[Theorem 1.1]{MR1718530} (see also \cite{MR1079188,MR1382207}). Because of the form of the functional $G_\rho,$ we cannot directly apply \cite{MR1718530}. But it can be easily adapted and we postpone its proof to the Appendix~\ref{appendixA} (see Theorem~\ref{thmJean} below). Let $I\subset(0,\infty)$ be any interval. Let us consider the following Property \eqref{PI}.
\begin{gather}
\label{PI}
\tag{$P_I$}
\begin{cases}
\text{For almost every } \rho\in I, \text{ there exists a sequence }	\medskip \\
(u_k)_{k\in\N} \subset H^1_{A,V}(\R^N) \text{ satisfying \eqref{eqPSB}--\eqref{eqPS}.}
\end{cases}
\end{gather}
Let $(a_n)_{n\in\N}\subset I_\gamma$ be any increasing sequence converging towards $\sup I_\gamma.$ For each $n\in\N,$ let $I_n\eqdef\left(\frac1n,a_n\right).$ Let $n\in\N$ be such that $I_n\neq\emptyset.$ By \eqref{MPG} and Theorem~\ref{thmJean}, for each $\rho_0\in I_\gamma,$ $I_{\rho_0}$ satisfies $(P_{I_{\rho_0}}).$ But $\ovl{I_n}\subset\bigcup\limits_{\rho_0\in I_\gamma}I_{\rho_0}$ and by compactness, $I_n$ may be covered by a finite number of $I_{\rho_0}.$ Consequently, $I_n$ satisfies $(P_{I_n}).$ Since $n\in\N,$ is  arbitrary, we infer that $I_\gamma=\bigcup\limits_{n\in\N}I_n$ satisfies $(P_{I_\gamma}).$ This ends the proof of the lemma.
\medskip
\end{proof*}

\begin{cor}
\label{corlemPS}
For almost every $\rho\in I_\gamma,$ there exists $u_\rho \in H^1_{A,V}(\R^N) \setminus\{0\}$ such that $G_\rho^\p (u_\rho) = 0.$ In particular, $u_\rho$ is a non zero weak solution to~\eqref{nlsp} with $\lambda=\frac1\rho.$
\end{cor}

\begin{proof*}
By Lemma~\ref{lemPS}, for almost every $\rho\in I_\gamma,$ there exist $c(\rho)>0$ and a sequence $(u_k)_{k\in\N}\subset H^1_{A,V}(\R^N)$ satisfying \eqref{eqPSB}--\eqref{eqPS}. Let such $\rho,$ $c\eqdef c(\rho)$ and $(u_k)_{k\in\N}.$ We first extract a subsequence (without change of notation) for which Theorem~\ref{thmTin} applies. By \eqref{eqPS}, the sequence $u_k\cancel{\xrightarrow[k\to\infty]{H^1_{A,V}}}0$ because $c>0$ and $G_\rho(0) =0.$ Thus we may assume that, up to a subsequence that we still denote by $(u_k)_{k\in\N},$ $\|u_k\|_{H^1_{A,V}(\R^N)}^2 \tends t > 0.$ It follows from \eqref{eqPS} that $\langle G_\rho^\p(u_k),u_k \rangle_{H^{-1}_{A,V},H^1_{A,V}} \tends 0.$ If $u_k\overset{D}{-\!\!\!-\!\!\!\weak}0$ then \eqref{thmLA2}, Proposition~\ref{propg}, \eqref{f1}--\eqref{f2}, Hölder's inequality, Sobolev's embedding and Theorem~\ref{thmcoc} imply that for any $\eps>0,$ there is a $p_\eps\in(2,\2)$ such that for some $C_\eps>0,$
\begin{align*}
	&	\; \left|\langle g(u_k),u_k\rangle_{H^{-1}_{A,V},H^1_{A,V}}\right|\le\vint_{\R^N}|g(u_k)||u_k|\d x				\\
  \le	&	\; C\left(\sup_{k\in\N}\|u_k\|_{H^1_{A,V}(\R^N)}\right)\eps+C_\eps\|u_k\|_{L^{p_\eps}(\R^N)}^{p_\eps}\xrightarrow{k\to\infty}C\eps.
\end{align*}
But then, $\langle G_\rho^\p(u_k),u_k\rangle_{H^{-1}_{A,V},H^1_{A,V}}\tends \rho t \neq 0,$ a contradiction. Then,
\begin{gather}
\label{demlemBPSsol}
u_k\cancel{\overset{D}{-\!\!\!-\!\!\!\weak}}0, \; \text{ as } k\tends\infty.
\end{gather}
Let us apply and use the notations of Theorem~\ref{thmTin}. If $\D=\emptyset$ or if all the $w^n$ were zero, then by \eqref{D_weak_limit} we would have $u_k\overset{D}{-\!\!\!-\!\!\!\weak}0,$ contradicting \eqref{demlemBPSsol}. Therefore, $\D\neq\emptyset$ and there is at least one nonzero $w^{n_0}$ which we call $u_\rho.$ Since for any $(k,n)\in\N\times\D,$
\begin{gather*}
G_\rho\left(\tau_{-y_k^n}u_k\right)=G_\rho(u_k) \; \text{ and } \;
\left\|G_\rho^\p\left(\tau_{-y_k^n}u_k\right)\right\|_{H^{-1}_{A,V}(\R^N)}=\|G_\rho^\p(u_k)\|_{H^{-1}_{A,V}(\R^N)},
\end{gather*}
we conclude from \eqref{eqPS}, \eqref{weak_limit}, Theorem~\ref{thmHHV} and Proposition~\ref{propg} that,
\begin{gather*}
G_\rho^\p\left(\tau_{-y_k^{n_0}}u_k\right)\underset{k\to\infty}{-\!\!\!-\!\!\!-\!\!\!\weak}G_\rho^\p(u_\rho)=0, \; \text{ in  }H^{-1}_{A,V}\text{-weakly},
\end{gather*}
from which the result follows.
\medskip
\end{proof*}

\begin{vproof}{of Theorem~\ref{thmmain}.}
Apply Corollary~\ref{corlemPS} with $\rho=\frac1\lambda.$
\medskip
\end{vproof}

\section{Applications}
\label{secapp}

In this section, we give some examples of nonlinearities for which Corollary~\ref{corlemPS} applies: for almost every $\lambda>0$ such that $\frac1\lambda\in I_\gamma=(0,2S),$ where $S\eqdef\vsup_{t \neq s} \frac{\gamma(t)-\gamma(s)}{t-s}\in(0,\infty],$ there exists, at least, a non zero weak solution to~\eqref{nlsp}.

\begin{exa}[\textbf{The single power interaction}]
\label{single}
Let $1<p<\2-1$ and let,
\begin{gather*}
\forall u\in H^1(\R^N), \; g(u)=|u|^{p-1}u.
\end{gather*}
Then \eqref{f1}--\eqref{psi+} are satisfied and Corollary~\ref{corlemPS} applies. It is not hard to see that $\frac{\gamma(t)}t\xrightarrow{t\to\infty}\infty$ so that $I_\gamma=(0,\infty).$ Let $\lambda>0.$ Let then $\lambda_0\in I_\gamma$ for which \eqref{nlsp} admits a non zero weak solution $u_{\lambda_0}.$ Setting $u=\left(\frac\lambda{\lambda_0}\right)^\frac1{p-1}u_{\lambda_0},$ a straightforward calculation shows that $u$ is a solution to~\eqref{nlsp} with $\lambda g(u)$ as the right side. In conclusion, for any $\lambda>0,$ equation~\eqref{nlsp} has, at least, a non zero weak solution. Note that $F$ satisfies the Rabinowitz condition.
\medskip
\end{exa}

\begin{exa}[\textbf{The combined power-type interaction}]
\label{double}
Let $\mu_1,\mu_2>0,$ let $1<p_1\neq p_2<\2-1$ and let,
\begin{gather*}
\forall u\in H^1(\R^N), \; g(u)=\mu_1|u|^{p_1-1}u-\mu_2|u|^{p_2-1}u.
\end{gather*}
The only difficulty is to show that there is a $u\in H^1(\R^N)$ such that,
\begin{gather*}
\psi(u)\eqdef\frac{\mu_1}{p_1+1}\|u\|_{L^{p_1+1}(\R^N)}^{p_1+1}-\frac{\mu_2}{p_2+1}\|u\|_{L^{p_2+1}(\R^N)}^{p_2+1}>0.
\end{gather*}
Let $u\in H^1_{A,V}(\R^N)\setminus\{0\}$ and let $t>0.$ If $p_1<p_2$ then
\begin{gather*}
\psi(tu)=t^{p_1+1}\left(\frac{\mu_1}{p_1+1}\|u\|_{L^{p_1+1}(\R^N)}^{p_1+1}-\frac{\mu_2}{p_2+1}t^{p_2-p_1}\|u\|_{L^{p_2+1}(\R^N)}^{p_2+1}\right)>0,
\end{gather*}
for any $0<t\ll1,$ while if $p_1>p_2$ then
\begin{gather*}
\psi(tu)=t^{p_2+1}\left(\frac{\mu_1}{p_1+1}t^{p_1-p_2}\|u\|_{L^{p_1+1}(\R^N)}^{p_1+1}-\frac{\mu_2}{p_2+1}\|u\|_{L^{p_2+1}(\R^N)}^{p_2+1}\right)
\xrightarrow{t\to\infty}\infty.
\end{gather*}
Then \eqref{f1}--\eqref{psi+} are satisfied and Corollary~\ref{corlemPS} applies. In particular, it follows from the last estimate that if $p_1>p_2$ then $\frac{\gamma(t)}t\xrightarrow{t\to\infty}\infty$ so that $I_\gamma=(0,\infty)$ and we may choose $\lambda$ as close to $1$ as we want. Notice also that $\psi<0$ on a nonempty open subset which is very different from the most hypotheses that can be found in the literature (as the Rabinowitz condition, for instance).
\medskip
\end{exa}

\begin{exa}
Suppose that \(F\) does not satisfy the Rabinowitz condition: \(F(x,t) \ge \mu tf(x,t)>0\) with \(\mu >2\) but there are an \(M >0\) and a $c>0$ such that for any \(t>M\), \(F(x,t) \ge ct^2 \ln t\). Then
\begin{align*}
  \lim_{t \to \infty} \frac{\gamma(t)}{t} & \ge c\lim_{t \to \infty}
                                               \sup_{u \in\vsS_t} \frac1t
                                               \int_{\R^N} |u|^2 \ln |u| \d x	\\
  & = c\lim_{t \to \infty} \sup_{u \in\vsS_1} \frac1t \int_{\R^N}t|u|^2 \ln(\sqrt{t}|u|) \d x \\
  & = \infty.
\end{align*}
So that \(I_\gamma = (0,\infty)\) and Corollary~\ref{corlemPS} applies for almost every $\lambda>0.$ As an example of $g$ satisfying such a condition and \eqref{f1}--\eqref{psi+} is,
\begin{gather*}
\forall u\in H^1(\R^N), \; g(u)=\
\begin{cases}
c_\eps|u|^{p-1}u,		&	\text{if } |u|<\eps,	\medskip \\
\mu_1u\ln|u|+\mu_2 u,	&	\text{if } |u|\ge\eps,
\end{cases}
\end{gather*}
where $\mu_1,\mu_2,\eps>0$ and $p\in(1,\2-1)$ can be chosen arbitrarily and $c_\eps=\eps^{-(p-1)}(\mu_1\ln\eps+\mu_2).$
\medskip
\end{exa}

\section*{Appendix}

\appendix
\section{Some proofs}
\label{appendixA}

In this appendix, we adapt the proof of \cite[Theorem 1.1]{MR1718530} to our family of functionals $(G_\rho)_{\rho\in I_\gamma},$ where the original idea is due to \cite{MR2431434}. We also give the proof of Proposition~\ref{propg}.
\medskip \\
In \cite{MR1718530}, the family of functionals is of the form
\begin{gather*}
\forall\lambda>0, \; I_\lambda(u)=A(u)-\lambda B(u),
\end{gather*}
where $A(u)\xrightarrow{\|u\|\to\infty}\infty$ or $B(u)\xrightarrow{\|u\|\to\infty}\infty,$ and with $B\ge0$ everywhere. Unfortunately, in our case,
\begin{gather*}
\forall\rho>0, \; G_\rho(u)=\frac1\lambda I_\lambda(u)=\frac1\lambda\left(\frac{\|u\|^2}2-\lambda\psi(u)\right), \;\; \lambda=\frac1\rho,
\end{gather*}
and we do not have $B=\psi\ge0,$ everywhere, but only somewhere. So we have, in some sense, to reverse the role of $A=\|u\|^2$ and $B=\psi.$ The following theorem is an easy adaptation of \cite[Theorem 1.1]{MR1718530}, but for the convenience of the reader, we give its proof.

\begin{thm}[\textbf{\cite{MR1718530}, Theorem 1.1}]
\label{thmJean}
Let $(X,\|\:.\:\|)$ be a Banach space, let $I\subset(0,\infty)$ be a nonempty open interval and let $(G_\rho)_{\rho\in I}\subset C^1\big(X;\R\big)$ be a family of functionals of the form,
\begin{gather}
\label{thmJean1}
\forall\rho\in I, \; G_\rho(u)=\rho A(u)-B(u),
\end{gather}
where $A\not\equiv0$ and for any $u\in X,$ $A(u)\ge0.$ Assume that either $A(u)\xrightarrow{\|u\|\to\infty}\infty$ or $B(u)\xrightarrow{\|u\|\to\infty}\infty.$ We also assume that $(G_\rho)_{\rho\in I}$ has mountain pass geometry: there exist $u_1\in X$ and $u_2\in X$ such that, denoting by 
\begin{gather*}
\Gamma\eqdef\left\{\xi\in C\left([0,1];X\right); \; \xi(0)=u_1 \text{ and } \xi(1)=u_2\right\},
\end{gather*}
the set of continuous paths joining $u_1$ to $u_2,$ we have for any $\rho\in I,$
\begin{gather}
\label{thmJean2}
c(\rho)\eqdef\inf_{\xi\in\Gamma}\max_{t\in[0,1]}G_\rho\big(\xi(t)\big)>\max\big\{G_\rho(u_1),G_\rho(u_2)\big\}.
\end{gather}
Then for almost every $\rho\in I,$ $G_\rho$ admits a bounded Palais-Smale sequence: there exists a sequence $(u_n)_{n\in\N}\subset X$ satisfying,
\begin{gather}
\label{thmJean3}
(u_n)_{n\in\N}\subset X \text{ is bounded,}			\\
\label{thmJean4}
\begin{cases}
G_\rho (u_n) \xrightarrow[n\to\infty]{} c(\rho),		\medskip \\
G_\rho^\p(u_n) \xrightarrow[n\to\infty]{X^\star}0,
\end{cases}
\end{gather}
where $X^\star$ denotes the topological space of $X.$
\end{thm}

\begin{rmk}
Here are some comments of Theorem~\ref{thmJean}.
\begin{enumerate}[1)]
\item
\label{rmkthmJean1}
If there exist $\rho\in I$ and $(u_1,u_2)\in X\times X$ satisfying \eqref{thmJean2} then it is well-known, by the Mountain Pass Theorem, that there exists a Palais-Smale sequence $(u_n)_{n\in\N}\subset X$ satisfying \eqref{thmJean4} (see, for instance, \cite[Theorem 6.2, p.144]{MR2294665}). The difficulty is to find such a \textbf{bounded} sequence.
\item
\label{rmkthmJean2}
The proof of Theorem~\ref{thmJean} relies on the existence of the derivative $c^\p(\rho)$ of $c(\rho).$ Since $A\ge0,$ we have by \eqref{thmJean2} that the mapping $c:\rho\longmapsto c(\rho)$ is nondecreasing over $I.$ It follows that $c$ has a derivative $c^\p$ almost everywhere on $I.$ In the original proof, the existence almost everywhere on $I$ of $c^\p$ is ensured by the fact that the mapping $c:\rho\longmapsto c(\rho)$ is nonincreasing over $I.$
\end{enumerate}
\medskip
\end{rmk}

\noindent
Before proceeding to the proof of Theorem~\ref{thmJean}, let us pick any $\rho\in I$ such that the derivative $c^\p(\rho)$ exists (see the item  \ref{rmkthmJean2}) in the above remark). Let then $\rho_0\in(0,\rho)$ be small enough to have $(\rho-\rho_0,\rho+\rho_0)\subset I$ and
\begin{gather}
\label{demthmJean0}
\forall\wt\rho\in(\rho-\rho_0,\rho+\rho_0), \; \left|\frac{c(\wt\rho)-c(\rho)}{\wt\rho-\rho}-c^\p(\rho)\right|\le1.
\end{gather}
Now, let us choose $(\rho_n)_{n\in\N}\subset(\rho,\rho+\rho_0)$ be a decreasing sequence such that $\rho_n\xrightarrow{n\to\infty}\rho.$ Finally, since $A(u)\xrightarrow{\|u\|\to\infty}\infty$ or $B(u)\xrightarrow{\|u\|\to\infty}\infty$ there exists $M>10$ such that for any $u\in X,$
\begin{gather}
\label{demthmJean00}
\|u\|>M \implies \max\big\{A(u),B(u)\big\}>\max\big\{c^\p(\rho)+3,2\rho\big(c^\p(\rho)+4\big)-c(\rho)\big\}.
\end{gather}
We shall need of the two following lemmas.

\begin{lem}
\label{lemJean1}
There exists $(\xi_n)_{n\in\N}\subset\Gamma$ satisfying the following properties.
\begin{enumerate}[$1)$]
\item
\label{lemJean11}
Let $t\in[0,1].$ If $n\in\N$ is such that $G_\rho\big(\xi_n(t)\big)\ge c(\rho)-(\rho_n-\rho)$ then $\|\xi_n(t)\|\le M.$
\item
\label{lemJean12}
$\forall n\in\N,$ $\vmax_{t\in[0,1]}G_\rho\big(\xi_n(t)\big)\le c(\rho)+(c^\p(\rho)+2)(\rho_n-\rho).$
\end{enumerate}
\end{lem}

\begin{proof*}
Let $(\xi_n)_{n\in\N}\subset\Gamma$ be such that for any $n\in\N,$
\begin{gather}
\label{demlemJean11}
\vmax_{t\in[0,1]}G_{\rho_n}\big(\xi_n(t)\big)\le c(\rho_n)+(\rho_n-\rho).
\end{gather}
Let $t\in[0,1].$ Let $n\in\N.$ We have by the hypothesis in \ref{lemJean11}), \eqref{demlemJean11} and \eqref{demthmJean0},
\begin{gather}
\label{demlemJean12}
A\big(\xi_n(t)\big)=\frac{G_{\rho_n}\big(\xi_n(t)\big)-G_\rho\big(\xi_n(t)\big)}{\rho_n-\rho}
\le\frac{c(\rho_n)-c(\rho)}{\rho_n-\rho}+2\le c^\p(\rho)+3.
\end{gather}
In addition, since for any $u\in X,$ the mapping $\rho\longmapsto G_\rho(u)$ is nondecreasing, it follows from \eqref{demlemJean12} and the hypothesis in \ref{lemJean11}),
\begin{gather}
\label{demlemJean13}
B\big(\xi_n(t)\big)=\rho_nA\big(\xi_n(t)\big)-G_{\rho_n}\big(\xi_n(t)\big)\le2\rho\big(c^\p(\rho)+4\big)-c(\rho).
\end{gather}
Hence $\|\xi_n(t)\|\le M,$ by \eqref{demthmJean00}, \eqref{demlemJean12} and \eqref{demlemJean13}. To prove the second part of the lemma, we see that \eqref{demthmJean0} implies,
\begin{gather}
\label{demlemJean14}
c(\rho_n)\le c(\rho)+\big(c^\p(\rho)+1\big)(\rho_n-\rho).
\end{gather}
Finally, \eqref{demlemJean11} and \eqref{demlemJean14} yield,
\begin{gather*}
\max_{t\in[0,1]}G_\rho\big(\xi_n(t)\big)\le\max_{t\in[0,1]}G_{\rho_n}\big(\xi_n(t)\big)\le c(\rho)+\big(c^\p(\rho)+2\big)(\rho_n-\rho).
\end{gather*}
This ends the proof of the lemma.
\medskip
\end{proof*}

\begin{lem}
\label{lemJean2}
Define for any $\eps>0,$
\begin{gather*}
F_\eps\eqdef\Big\{u\in X; \; \|u\|\le 2M \text{ and } |G_\rho(u)-c(\rho)|\le\eps\Big\}.
\end{gather*}
Then for any $\eps>0,$ $F_\eps\neq\emptyset$ and $\vinf_{u\in F_\eps}\|G_\rho^\p(u)\|_{X^\star}=0.$
\end{lem}

\begin{proof*}
Let $(\xi_n)_{n\in\N}\subset\Gamma$ be given by Lemma~\ref{lemJean1}. Then for each $n\in\N,$ there exists $t_n\in[0,1]$ such that $0\le G_\rho\big(\xi_n(t_n)\big)-c(\rho)\le(c^\p(\rho)+2)(\rho_n-\rho)\xrightarrow{n\to\infty}0$ and $\|\xi(t_n)\|\le M.$ We infer that for any $\eps>0,$ there exists $n_0\in\N$ large enough such that $\xi(t_{n_0})\in F_\eps.$ Now, we note that it is sufficient to show the result for any $\eps>0$ small enough. If the result does not hold then there exists $0<\eps_0<\frac{c(\rho)-\max\{G_\rho(u_1),G_\rho(u_2)\}}2$ such that $\vinf_{u\in F_{2\eps_0}}\|G_\rho^\p(u)\|_{X^\star}\ge2\eps_0.$ We then may apply a deformation lemma to affirm that there exists a homeomorphism $\eta:X\tends X$ satisfying the following properties.
\begin{align}
\label{demlemJean21}
&	\text{If } |G_\rho(u)-c(\rho)|>2\eps_0 \text{ then } \eta(u)=u.	\\
\label{demlemJean22}
&	\forall u\in X, \; G_\rho\big(\eta(u)\big)\le G_\rho(u).			\\
\label{demlemJean23}
&	\text{If } \|u\|\le M \text{ and } G_\rho(u)<c(\rho)+\eps_0 \text{ then } G_\rho\big(\eta(u)\big)<c(\rho)-\eps_0.
\end{align}
See for instance \cite[Theorem 4.2, p.38]{MR2012778}. The assertion \eqref{demlemJean22} is not directly stated in this theorem but in its proof p.39. Let $m\in\N$ be large enough to have,
\begin{gather}
\label{demlemJean24}
\rho_m-\rho<(c^\p(\rho)+2)(\rho_m-\rho)<\eps_0.
\end{gather}
By \eqref{demlemJean21}, $\eta(\xi_m)\in\Gamma.$ Let $t\in[0,1].$

$\bullet$
If $G_\rho\big(\xi_m(t)\big)\le c(\rho)-(\rho_m-\rho)$ then by \eqref{demlemJean22},
\begin{gather}
\label{demlemJean25}
G_\rho\big(\eta(\xi_m(t))\big)\le c(\rho)-(\rho_m-\rho).
\end{gather}

$\bullet$
If $G_\rho\big(\xi_m(t)\big)>c(\rho)-(\rho_m-\rho)$ then by Lemma~\ref{lemJean1} and \eqref{demlemJean24}, $\|\xi_m(t)\|\le M$ and $G_\rho\big(\xi_m(t)\big)<c(\rho)+\eps_0.$ It then follows from \eqref{demlemJean23} and \eqref{demlemJean24},
\begin{gather}
\label{demlemJean26}
G_\rho\big(\eta(\xi_m(t))\big)<c(\rho)-\eps_0<c(\rho)-(\rho_m-\rho).
\end{gather}
It follows from \eqref{demlemJean25} and \eqref{demlemJean26} that,
\begin{gather*}
c(\rho)=\inf_{\xi\in\Gamma}\max_{t\in[0,1]}G_\rho\big(\xi(t)\big)\le\max_{t\in[0,1]}G_\rho\big(\eta(\xi_m(t))\big)\le c(\rho)-(\rho_m-\rho).
\end{gather*}
A contradiction, since $\rho_m-\rho>0.$
\medskip
\end{proof*}

\begin{vproof}{of Theorem~\ref{thmJean}.}
The result follows by applying Lemma~\ref{lemJean2} with any sequence $\eps_n\searrow0.$
\medskip
\end{vproof}

\begin{vproof}{of Proposition~\ref{propg}.}
Throughout this proof, we let $\kappa=1,$ if $N\ge3$ and $\kappa=0,$ if $N\le2.$ We will denote by $C_1>1$ and $p_1$ the constants given by \eqref{f1}--\eqref{f2} for $\eps=1.$ We proceed to the proof in 6 steps.
\\
\textbf{Step 1:} $g:H^1(\R^N)\tends H^{-1}(\R^N)$ is well-defined, bounded on bounded sets and \ref{propg2} holds.
\\
By \eqref{f1}--\eqref{f2}, $g(u)\in L^1_\loc(\R^N).$ Let $\vphi\in\Dr(\R^N).$ We have by \eqref{f1}--\eqref{f2}, Hölder's inequality and the Sobolev embeddings,
\begin{align*}
	&	\; \left|\langle g(u),\vphi\rangle_{\Dr^\p(\R^N);\Dr(\R^N)}\right|=\left|\Re\dsp\vint_{\R^N}g(u)\ovl\vphi\d x\right|					\\
  \le	&	\; C_1\left(\|u\|_{L^2(\R^N)}+\kappa\|u\|_{L^\2(\R^N)}^{\2-1}+\|u\|_{L^{p_1}(\R^N)}^{p_1-1}\right)\|\vphi\|_{H^1(\R^N)}		\\
  \le	&	\; C\left(\|u\|_{H^1(\R^N)}+\kappa\|u\|_{H^1(\R^N)}^{\2-1}+\|u\|_{H^1(\R^N)}^{p_1-1}\right)\|\vphi\|_{H^1(\R^N)}.
\end{align*}
By density, it follows that $g:H^1(\R^N)\tends H^{-1}(\R^N)$ is well-defined, $g$ is bounded on bounded sets and Property~\ref{propg2} holds. \\
\textbf{Step 2:} $\psi\in C(H^1(\R^N);\R),$ $\psi$ is bounded on bounded sets, Gâteaux-differentiable and its Gâteaux-differential is $\psi^\p_\g=g.$ \\
Let $u\in H^1(\R^N).$ By \eqref{f1}--\eqref{f2}, Hölder's inequality and the Sobolev embedding, $F(u)\in L^1(\R^N;\R)$ so that $\psi:H^1(\R^N)\tends\R$ is well-defined and $\psi$ is bounded on bounded sets. Let $v\in H^1(\R^N).$ Still by \eqref{f1}--\eqref{f2}, Hölder's inequality and the Sobolev embedding,
\begin{align*}
	&	\; |\psi(u+v)-\psi(u)|\le\vint_{\R^N}\vint_{|u|}^{|u+v|}\big(t+\kappa t^{\2-1}+C_1t^{p_1-1}\big)\d t\d x	\\
  \le	&	\; C\Big(\|u\|_{L^2}+\|v\|_{L^2}+\kappa(\|u\|_{L^\2}+\|v\|_{L^\2})^{\2-1}+(\|u\|_{L^{p_1}}
										+\|v\|_{L^{p_1}})^{p_1-1}\Big)\|v\|_{H^1(\R^N)}.
\end{align*}
It follows that $\psi\in C(H^1(\R^N);\R).$ Let $v\in H^1(\R^N)$ and $0<|t|<1.$ Since $u,v\in L^2(\R^N),$ the set
\begin{gather*}
\vN\eqdef\big\{x\in\R^N; |u(x)|=\infty \text{ or } |v(x)|=\infty\big\},
\end{gather*}
has Lebesgue measure 0. Let $x\in\vN^\co.$ If $u(x)\neq0$ then using that
\begin{gather*}
|u(x)+tv(x)|=\sqrt{\big(u(x)+tv(x)\big)\ovl{\big(u(x)+tv(x)\big)}}>0,
\end{gather*}
for  $t$ small enough, we  see that
\begin{gather*}
\frac\d{\d t}F(x,|u(x)+tv(x)|)_{|t=0}=\Re\big(\f(x,u(x))\ovl{v(x)}\big).
\end{gather*}
If $u(x)=0$ then by \eqref{f1}--\eqref{f2},
\begin{gather*}
\left|\frac{F(x,|tv(x)|)-F(x,0)}t\right|\le C\big(|t||v(x)|^2+\kappa|t|^{\2-1}|v(x)|^\2+C_1|t|^{p_1-1}|v(x)|^{p_1}\big)\xrightarrow{t\to0}0.
\end{gather*}
We then infer,
\begin{gather*}
\frac{F(\:.\:,|u+tv|)-F(\:.\:,|u|)}t\xrightarrow[t\tends0]{\text{a.e. in }\R^N}\Re\big(\f(\:.\:,u)\ovl v\big).
\end{gather*}
By \eqref{f1}--\eqref{f2},
\begin{align*}
	&	\; \frac{F(\:.\:,|u+tv|)-F(\:.\:,|u|)-t\Re\big(\f(x,u)\ovl v\big)}t						\\
  \le	&	\; \frac1t\vint_{|u|}^{|u+tv|}|f(\:.\:,s)|\d s+|f(\:.\:,|u|)||v|							\\
  \le	&	\; C\big(|u|+|v|+\kappa(|u|+|v|)^{\2-1}+(|u|+|v|)^{p_1-1}\big)|v|\in L^1(\R^N).
\end{align*}
It follows from the dominated convergence Theorem and Property~\ref{propg2} that,
\begin{gather*}
\lim_{t\to0}\frac{\psi(u+tv)-\psi(u)}t=\langle g(u),v\rangle_{H^{-1}(\R^N),H^1(\R^N)}.
\end{gather*}
Hence Step~2.
\\
\textbf{Step 3:} Let $u,v\in H^1(\R^N)$ and $(u_n)_{n\in\N}\subset H^1(\R^N)$ be bounded. Let $\eps>0.$ Choose $\eps^\p>0$ small enough to have,
\begin{gather}
\label{dempropg1}
2\eps^\p\left(\sup_{n\in\N}\|u_n\|_{L^2(\R^N)}^2+\|u\|_{L^2(\R^N)}^2+\kappa\left(\sup_{n\in\N}\|u_n\|_{L^\2(\R^N)}^{\2-1}+\|u\|_{L^\2(\R^N)}^{\2-1}\right)\right)\le\eps.
\end{gather}
For such an $\eps^\p,$ let $p_{\eps^\p}$ and $C_{\eps^\p}$ be given by \eqref{f1}--\eqref{f2}. For each $n\in\N,$ let
\begin{gather*}
A_n=\Big\{x\in\R^N;\eps^\p\big(|u_n|+|u|+\kappa(|u_n|^{\2-1}+|u|^{\2-1})\big)\le C_{\eps^\p}(|u_n|^{p_{\eps^\p}-1}+|u|^{p_{\eps^\p}-1})\Big\}.
\end{gather*}
It holds that,
\begin{gather}
\label{dempropg2}
\left|\langle g(u_n)-g(u),v\rangle_{H^{-1}(\R^N),H^1(\R^N)}\right|
\le\vint_{\R^N}\big|g(u_n)-g(u)\big||v|\1_{A_n}\d x+\eps\|v\|_{H^1(\R^N)}.
\end{gather}
Indeed, by \eqref{f1}--\eqref{f2}, Hölder's inequality, the Sobolev embeddings and \eqref{dempropg1}, we have,
\begin{align*}
	&	\; \left|\langle g(u_n)-g(u),v\rangle_{H^{-1}(\R^N),H^1(\R^N)}\right|\le\vint_{\R^N}\big|g\big(u_n)-g(u)\big||v|\d x		\\
   =	&	\; \vint_{\R^N}\big|g(u_n)-g(u)\big||v|\1_{A_n}\d x+\vint_{\R^N}\big|g\big(u_n)-g(u)\big||v|\1_{A_n^\co}\d x			\\
  \le	&	\; \vint_{\R^N}\big|g(u_n)-g(u)\big||v|\1_{A_n}\d x
			+2\eps^\p\vint_{\R^N}\big(|u_n|+|u|+\kappa(|u_n|^{\2-1}+|u|^{\2-1})\big)|v|\d x							\\
  \le	&	\; \vint_{\R^N}\big|g(u_n)-g(u)\big||v|\1_{A_n}\d x+\eps\|v\|_{H^1(\R^N)}.
\end{align*}
Step~3 is proved.
\\
\textbf{Step~4:} $\psi\in C^1(H^1(\R^N);\R)$ and $\psi^\p=g.$
\\
By Step~2, it remains to show that $g\in C(H^1(\R^N);H^{-1}(\R^N))$ to have that $\psi$ is Fréchet-differentiable and $\psi^\p=\psi^\p_\g.$ Assume $u_n\xrightarrow[k\to\infty]{H^1(\R^N)}u.$ Let $\eps>0.$ Let then $\eps^\p,$ $p_{\eps^\p}$ and $C_{\eps^\p}$ be given by Step~3. By Hölder's inequality, we have for any $v\in H^1(\R^N),$
\begin{gather}
\label{dempropg3}
\vint_{\R^N}\big|g(u_n)-g(u)\big||v|\1_{A_n}\d x\le\big\|\big(g(u_n)-g(u)\big)\1_{A_n}\big\|_{L^{p^\p_{\eps^\p}}(\R^N)}\|v\|_{L^{p_{\eps^\p}}(\R^N)}.
\end{gather}
It follows from Sobolev' embedding and \eqref{dempropg2}--\eqref{dempropg3} that,
\begin{gather}
\label{dempropg4}
\sup_{\|v\|_{H^1(\R^N)}=1}\left|\langle g(u_n)-g(u),v\rangle_{H^{-1}(\R^N),H^1(\R^N)}\right|
\le C\big\|\big(g(u_n)-g(u)\big)\1_{A_n}\|_{L^{p^\p_{\eps^\p}}(\R^N)}+\eps
\end{gather}
We claim that,
\begin{gather}
\label{dempropg5}
\lim_{n\to\infty}\big\|\big(g(u_n)-g(u)\big)\1_{A_n}\|_{L^{p^\p_{\eps^\p}}(\R^N)}=0.
\end{gather}
If not, for some $\eps_0>0$ and a subsequence, that we will denote by $(u_n)_n,$ there would exist $h\in L^{p_{\eps^\p}}(\R^N;\R)$ such that for any $n\in\N,$ $\big\|\big(g(u_n)-g(u)\big)\1_{A_n}\|_{L^{p^\p_{\eps^\p}}(\R^N)}\ge\eps_0,$ $|u_n|\overset{\text{a.e}}{\le}h$ and $u_n\xrightarrow[n\to\infty]{\text{a.e. in }\R^N}u.$ But then $\big(g(u_n)-g(u)\big)\1_{A_n}\xrightarrow[n\tends\infty]{\text{a.e. in }\R^N}0$ and $\big|g(u_n)-g(u)\big|\1_{A_n}\le Ch^{p_{\eps^\p}-1}\in L^{p^\p_{\eps^\p}}(\R^N).$ This would yield to a contradiction by the Lebesgue convergence Theorem. Hence \eqref{dempropg5}. It then follows from \eqref{dempropg4}--\eqref{dempropg5} that,
\begin{gather*}
\forall\eps>0, \;
\limsup_{n\to\infty}\|g(u_n)-g(u)\|_{H^{-1}(\R^N)}\le\eps.
\end{gather*}
Letting $\eps\searrow0,$ we get $g\in C(H^1(\R^N);H^{-1}(\R^N)).$
\\
\textbf{Step 5}: Let $(u_n)_n,(v_n)_n \subset H^1(\R^N)$ be bounded. If $\vlim_{n\to\infty}\|u_n-v_n\|_{L^p(\R^N)}=0,$ for some $p\in[1,\infty],$ then $\vlim_{n\to\infty}|\psi(u_n)-\psi(v_n)|=0.$
\\
Let $\eps>0.$ For such an $\eps,$ let $p_\eps$ and $C_\eps$ be given by \eqref{f1}--\eqref{f2}. Let for any $t\in[0,1],$ $a(t)=\psi(v_n+t(u_n-v_n)).$ Then $a\in C^1([0,1];\R)$ and by the mean value Theorem, there exists $t_n\in(0,1)$ such that $a(1)-a(0)=a^\p(t_n)(1-0),$ that is
\begin{gather*}
\psi(u_n)-\psi(v_n)=\langle g(w_n),u_n-v_n\rangle_{H^{-1}(\R^N),H^1(\R^N)}.
\end{gather*}
where $w_n=v_n+t_n(u_n-v_n).$ Note that $(w_n)_{n\in\N}$ is bounded in $H^1(\R^N).$ It follows from \eqref{f1}--\eqref{f2}, Hölder's inequality and Sobolev's embedding that $\vlim_{n\to\infty}\|u_n-v_n\|_{L^{p_\eps}(\R^N)}=0$ and
\begin{align*}
	&	\; |\psi(u_n)-\psi(v_n)|																\\
  \le	&	\; \eps\left(\|w_n\|_{L^2(\R^N)}+\kappa\|w_n\|_{L^\2(\R^N)}^{\2-1}\right)\|u_n-v_n\|_{H^1(\R^N)}
					+C_{p_\eps}\|w_n\|_{L^{p_\eps}(\R^N)}^{{p_\eps}-1}\|u_n-v_n\|_{L^{p_\eps}(\R^N)}	\\
  \le	&	\; C\eps+C\|u_n-v_n\|_{L^{p_\eps}(\R^N)}.
\end{align*}
We infer,
\begin{gather*}
\forall\eps>0, \; \limsup_{n\to\infty}|\psi(u_n)-\psi(v_n)|\le C\eps,
\end{gather*}
from which the result follows.
\\
\textbf{Step 6:} If $u_n \overset{H^1_\w}{-\!\!\!-\!\!\!-\!\!\!\weak}u$ then $g(u_n) \overset{H^{-1}_\w}{-\!\!\!-\!\!\!-\!\!\!\weak} g(u).$
\\
Since $(g(u_n))_{n\in\N}$ is bounded in $H^{-1}(\R^N)$ (Step~1), it is enough to show that $g(u_n)\xrightarrow[n\to\infty]{\Dr^\p(\R^N)}g(u).$ Let $\vphi\in\Dr(\R^N)$ with $\supp\vphi\subset B(0,R),$ for some $R>0.$ By compactness, $u_n\xrightarrow[n\to\infty]{L^{p_{\eps^\p}}(B(0,R))}u.$ Arguing by contradiction and using the dominated convergence Theorem, we show in the same way as in Step~4,
\begin{gather*}
\lim_{n\to\infty}\vint_{\R^N}\big|g(u_n)-g(u)\big||\vphi|\1_{A_n}\d x=0,
\end{gather*}
from which we deduce, with help of  \eqref{dempropg2},
\begin{gather*}
\forall\eps>0, \;
\limsup_{n\to\infty}\left|\langle g(u_n)-g(u),\vphi\rangle_{\Dr^\p(\R^N);\Dr(\R^N)}\right|\le\eps\|\vphi\|_{H^1(\R^N)}.
\end{gather*}
We conclude as in Step~4.
\medskip
\end{vproof}

\section{Topological vector spaces over the field of complex numbers restricted to the field of real numbers}
\label{appendixB}

Throughout this paper, we consider Banach spaces (or, more generally, complete topological vector spaces) over $\R$ rather than $\C.$ The main motivations are the following. Firstly, the linear forms are real-valued and there is a relation of order over $\R.$ Secondly, if a function $\psi$ belongs to $C^1(X;\R)$ (as in Proposition~\ref{propg}, for instance), where $X$ is a real Banach space, then $\psi^\p\in C(X;X^\star),$ where $X^\star$ is the $\R$-vector space $\Lr(X;\R).$ If $X$ is a complex Banach space then $X^\star$ is the $\C$-vector space $\Lr(X;\C)$ and $\psi^\p\in C\big(X;\Lr(X;\R)\big).$ But then, when a Riesz representation theorem exists, we have two kinds of representation between the elements of $\Lr(X;\R)$ and those of $X^\star=\Lr(X;\C),$ since $\Lr(X;\R)$ is not $\C$-linear. On the other hand, if $X$ is a complex Banach space, it could be pleasant to consider $\lambda x,$ for $(\lambda,x)\in\C\times X.$ So, if $X_\C$ is a complex topological vector space, throughout this paper we consider $X_\R$ as the elements of $X_\C$ over the field $\R.$ We then consider the  real topological vector space $X^\star_\R$. For any $(\lambda,x)\in\C\times X,$ $\lambda x\in X_\R,$ since $X_\R$ and $X_\C$ have the same elements. In the special case where $H_\C$ is a complex Hilbert space whose the inner product is $(\:.\:,\:.\:)_H$ then $H_\R$ is the real Hilbert space whose the scalar product is $\langle\:.\:,\:.\:\rangle_H\eqdef\Re\:(\:.\:,\:.\:)_H.$ In particular, for any $(u,v)\in H_\R\times H_\R,$ $\langle\vi u,\vi v\rangle_H=\langle u,v\rangle_H.$ Now, assume that $X_\C$ is a complex Banach space. Denote by $X^\star_\C$ and $X^\star_\R$ the topological dual spaces of $X_\C$ and $X_\R,$ respectively. It follows that $X^\star_\C$ is a $\C$-linear space while $X^\star_\R$ is only a $\R$-linear space. Let us define the map,
\begin{gather}
\label{I}
\begin{array}{rcl}
I:X^\star_\C	&	     \tends		&	X^\star_\R,	\medskip \\
		L	&	\longmapsto	&	\Re\:L.
\end{array}
\end{gather}
Then $I$ is a bijective isometry from $X^\star_\C$ onto $X^\star_\R$ (Brezis~\cite[Proposition~11.22, p.361]{MR2759829}). With help of this correspondance, we can identify some linear forms. For instance, let $X=L^p(\Omega;\C),$ where $\Omega$ is an open subset of $\R^N$ and $1\le p<\infty.$ If $p=2$ then the inner and scalar products are given by
\begin{gather*}
(u,v)_X=\vint_\Omega u(x)\ovl{v(x)}\d x \quad \text{ and } \quad \langle u,v\rangle_X=\Re\vint_\Omega u(x)\ovl{v(x)}\d x,
\end{gather*}
respectively. Using the  Riesz representation Theorem for the complex $L^p(\Omega;\C)_\C$ spaces (Yosida~\cite[Example~3, p.115]{MR617913}) and the bijective isometric map~\eqref{I}, it follows that
\begin{gather*}
L^p(\Omega;\C)^\star_\R=L^{p^\p}(\Omega;\C)_\R,
\end{gather*}
where $\dfrac1p+\dfrac1{p^\p}=1.$ More precisely, for any $L\in L^p(\Omega;\C)^\star_\R,$ there exists a unique $u\in L^{p^\p}(\Omega;\C)_\R$ such that
\begin{gather*}
\langle L,v\rangle_{L^p(\Omega)^\star,L^p(\Omega)}=\Re\vint_\Omega u(x)\ovl{v(x)}\d x,
\end{gather*}
for any $v\in L^p(\Omega;\C)_\R.$ Furthermore, $\|u\|_{L^{p^\p}(\Omega;\C)_\R}=\|L\|_{L^p(\Omega;\C)^\star_\R}.$ Finally, we end this appendix with the space of distributions $\Dr^\p(\Omega;\C).$ We consider the $\C$-complete topological vector space $\Dr(\Omega;\C)$ restricted to the field $\R$ as above. Then an element $T$ belongs to the $\R$-complete topological vector space $\Dr^\p(\Omega;\C)$ if $T$ is a $\R$-linear continuous mapping from $\Dr(\Omega;\C)$ to $\R.$ In particular, a function $f\in L^1_\loc(\Omega;\C)$ (over the field $\R)$ defines a distribution $T_f\in\Dr^\p(\Omega;\C)$ by the formula,
\begin{gather*}
\langle T_f,\vphi\rangle_{\Dr^\p(\Omega;\C),\Dr(\Omega;\C)}=\Re\vint_\Omega f(x)\ovl{\vphi(x)}\d x,
\end{gather*}
for any $\vphi\in\Dr(\Omega;\C).$ Indeed, $T_f$ is clearly a $\R$-linear continuous mapping from $\Dr(\Omega;\C)$ to $\R.$ Furthermore, if $f\in L^1_\loc(\Omega;\C)$ satisfies,
\begin{gather*}
\Re\vint_\Omega f(x)\ovl{\vphi(x)}\d x=0,
\end{gather*}
for any $\vphi\in\Dr(\Omega;\C),$ then $f=0.$ To see this, we note that $\Re(f),\Im(f)\in L^1_\loc(\Omega;\R)$ and choosing $\vphi=\psi+\vi0$ and then $\vphi=0+\vi\psi$ in the above expression, we get
\begin{gather*}
\vint_\Omega\Re\big(f(x)\big)\psi(x)\d x=\vint_\Omega\Im\big(f(x)\big)\psi(x)\d x=0,
\end{gather*}
for any $\psi\in\Dr(\Omega;\R).$ We infer that $\Re(f)=\Im(f)=0$ (Brezis~\cite[Corollary~4.24, p.110]{MR2759829}), from which the result follows. Obviously, if $f_n\xrightarrow[n\to\infty]{L^1_\loc(\Omega;\C)}f$ then $T_{f_n}\xrightarrow[n\to\infty]{\Dr^\p(\Omega;\C)}T_f.$ We conclude that,
\begin{gather*}
L^1_\loc(\Omega;\C)\inj \Dr^\p(\Omega;\C),
\end{gather*}
with embedding $T:f\in L^1_\loc(\Omega;\C)\longmapsto T_f\in\Dr^\p(\Omega;\C).$

\section*{Acknowledgement}
The authors would like to thank the anonymous referees for their valuable feedback and recommendations.

\end{document}